
\input amstex
\magnification=1200
\documentstyle{amsppt}
\NoBlackBoxes
\documentstyle{amsppt}
\centerline{ \bf Spinor Representations of $U_q(\hat{\frak gl}(n))$
and}
\centerline{\bf Quantum  Boson-Fermion Correspondence}
\vskip .2in
\centerline{Jintai Ding}
\centerline{RIMS, Kyoto University}
\vskip .4in
\centerline{ abstract}
This is an extension of quantum spinor construction in \cite{DF2}. We
define quantum affine Clifford algebras based on the tensor category and
the solutions of q-KZ equations, construct quantum spinor representations
of $U_q(\hat{\frak gl}(n))$ and explain classical and quantum
boson-fermion correspondence.

\subheading{I. Introduction}

The independent discovery of a  q-deformation
of universal enveloping algebra of an arbitrary Kac-Moody algebra
by Drinfeld \cite{D1} and Jimbo \cite{J1}
immediately raised numerous questions about
q-deformations of various structures associated to Kac-Moody algebras. A major
step in this direction was the work  by
Lusztig \cite{L}, who obtained a q-deformation of the
category of highest weight representations of Kac-Moody algebras for
generic or formal parameter q.
There were a number of
successful results on q-deformation of various
mathematical structures of finite dimensional and
affine Lie algebras  \cite{D2},
 \cite{J2}, \cite{FJ}, \cite{H}, \cite{FR}, etc. However
each particular problem required its own special insight
 and in some cases presented formidable
difficulties.

In \cite{DF2}, we propose an invariant approach,
which was also stressed in  \cite{FRT}, where a q-analogue of matrix
realization of classical Lie algebras was given. We manage to use such
an  approach to
define quantum Clifford and Weyl algebras
using general representation theory of quantum groups.
We  show that the explicit formulas for quantum Clifford and Weyl algebras
 match the ones actively studied in physics literature
(see e.g. \cite{WZ}, \cite{K}). Using
those quantum algebras, we construct spinor and oscillator representations of
quantum groups of classical types and recover all the relevant formulas for the
quantum construction.  Uniqueness
arguments from representation theory allow us to justify
that the quadratic
expressions in quantum Clifford and Weyl algebras provide the
desired representations. An explicit verification
of Serre's relations for quantum groups lead to rather
involved formulas (cf. \cite{H}).
 The key idea   consists of reformulating
 familiar classical constructions entirely
in terms of the  tensor category of highest weight representations and
then,
using Lusztig's result on q-deformation of this category to  define the
corresponding quantum structures. In the quantum case, we often need a
quasitriangular structure of the  tensor category introduced by Drinfeld
\cite{D3},
since it plays the role of the symmetric structure in the classical
case. In particular, we would like to emphasize
the central role of the universal Casimir operator implied by
the quantum structure.

The motivation to develop  such an invariant approach for the
quantum groups corresponding to the simple finite dimensional algebras
is surely to apply it to general cases, especially
 affine quantum groups. In this paper, we  parallelly extend
 the idea  to the cases of  the
spinor representations of
quantum affine groups $U_q(\hat {\frak gl}(n))$
and $U_q(\hat {\frak sl}(n))$,  though there is  slight  difference.
In a  subsequent work, we will present  the
corresponding results of the cases of $\hat {\frak o}(n)$ and
$\hat{\frak sp}(2n)$ as for the cases of ${\frak o}(n)$ and
${\frak sp}(2n)$ in \cite{DF2}.

For undeformed  affine Lie algebras,
Garland showed  \cite{G} that the affine
Kac-Moody algebra $\hat {\frak{g}}$ associated to a simple Lie algebra
$\frak{g}$ admits a natural realization as a central extension of the
corresponding loop algebra ${\frak{ g}} \otimes \Bbb C [t, t^{-1}] $.

For classical affine Lie algebras, we have the Frenkel-Kac construction to
obtain so-called  bosonic representations by vertex operators in terms of
Heisenberg algebra. There also exists another family of
representations,
the fermionic representations \cite{F} \cite{FF} in terms  of
affine Clifford algebras.  For these two families of representations,
there exists the well known  boson-fermion correspondence
\cite{F1} to show their equivalent relations.

As  in the undeformed classical cases, with the canonical basis for
the n-dimension  representation of ${\frak gl}(n)$, starting from
abstract representations without  any concrete realization,
 we can write down formulas
that  completely parallel
to the ones in \cite{DF2} to obtain the isomorphism between the
algebra generated by intertwiners and affine Clifford algebra,  and
  derive the  spinor construction. This
 provides a structural  explanation of the
boson-fermion correspondence.

For the quantum case, Faddeev, Reshetikhin and
Takhtajan \cite{FRT2} have shown how to extend their realization of
$U_q(\frak{g})$ to the quantum loop algebra $U_q(
{\frak{ g}} \otimes \Bbb [t, t^{-1}])$ via a canonical
 solution
 of the Yang-Baxter equation depending on a parameter $z\in\Bbb C$.
The first realization of the quantum affine algebra $U_q(\hat{
\frak g})$ and its special degeneration called the Yangian were
 obtained by Drinfeld \cite{D2}.

Later Reshetikhin and Semenov-Tian-Shansky
\cite{RS}  incorporated  the central extension in the
previous construction of \cite{FRT2} to obtain  the second
realization of the quantum affine algebra $U_q(\hat{ \frak{g}})$.

As  for the cases of deformation of finite dimensional Lie algebras,
with the theory of Lusztig of the q-deformation of the category of
highest weight representations,
we  construct the spinor representation of $U_q(\hat {\frak gl}(n))$
with  intertwiners completely
parallel to the case of the quantum groups of classical types as
explained in \cite{DF2}.

On the other hand, with the  Heisenberg algebra, Frenkel and Jing \cite{FJ}
constructed representations of $U_q(\hat {\frak sl}(n))$
by vertex operators in terms of Drinfeld's realization, which
could be extended to $U_q(\hat {\frak gl}(n))$.
This  gives us the quantized bosonic construction.

In \cite{DF}, we explicitly  established  a relation
between   two realizations of the quantum affine algebra $U_q(\hat
{\frak{gl}}(n))$. We show that the realization of Drinfeld's construction
can be naturally established in the Gauss decomposition of a matrix composed
of elements of the quantum affine algebra.
With this isomorphism, we are able to establish the  quantum boson-fermion
correspondence for the case of $U_q(\hat{\frak gl}(n))$
by the   concrete  isomorphism between those two representations
via the intertwiners, whose bosonization is partially solved in \cite{Ko}.

Hayashi constructed fermionic representations of $U_q(\hat{\frak
sl}(n))$, where however he only obtained the realization of the
basic generators. From the point view of this paper,
his  definition of  the
quantized fermions does not reveal the new structure
quantum fermions possess.

The quantum Clifford algebras we define are closely related to
massive quantum field theory\cite{Sm}. We will construct a realization
of an algebra, which resmebles the algebra to define
form factors in massive quantum field theory based on
the quantum Clifford algebra.

This paper  is arranged as follows. We present the basic idea of this paper
in Section 1. Then in  Section 2,  we will
present different realizations of $U_q(\hat{\frak gl}(n))$.
In section 3, we will reconstruct classical spinor representations
and boson-fermion correspondence. Section 4  will give the construction
of spinor  representations of $U_q(\hat {\frak gl}(n))$ and
explain the   quantum   boson-fermion correspondence.
We will also discuss the connection with the theory of
form factors in massive quantum field theory.

\documentstyle{amsppt}
{\bf 1.  Univeral Casimir Operators of
Affine Quantum Groups}

 The definitions  of quantum groups
corresponding to the affine Lie algebras, which were  given by Drinfeld
\cite{D1} and Jimbo \cite{J1} in terms of
generators and quantized relations corresponding to the affine Cartan
matrices, are very simple.

\proclaim{Definition 1.1}  Let $(a_{ij})$ be a  Cartan
matrix for an  affine Lie algebra $\hat {\frak g}$ \cite{Ka}. Let
 $(d_0,\ldots ,d_n)$ be a vector with integer entries $d_i\in
\{ 1,2,3,4\}$ such that $(d_i a_{ij})$ is symmetric.
 Let $q$ be an indeterminate.
  $U_q(\hat {\frak g})$ is an associative algebra on
$\Bbb C[q,q^{-1}]$
 with generators $E_i, F_i, K_i$,$K_i^{-1}$($i=0,1..,n$)
and the relations are:
$$\align
K_i K_j &= K_j K_i \ \
K_i K_i^{-1} = K_i^{-1} K_i = 1 \\
K_i E_j &= q^{d_ia_{ij}}E_j K_i, \qquad K_i F_j = q^{-d_ia_{ij}}
F_jK_i\\
E_i E_j &- F_j E_i = \delta_{ij}
\frac{K_{i}^{}-K_{j}^{-1}}{q^{d_i}-q^{-d_i}}
\endalign
$$
$$\align
&\sum\limits_{r+s=1-a_{ij}} (-1)^s \bmatrix 1-a_{ij} \\
s\endbmatrix_{d_{i}} E^r_i E_j E^s_i = 0,  \qquad if~~ i\neq j,  \tag 1.1\\
&\sum\limits_{r+s=1-a_{ij}} (-1)^s \bmatrix 1-a_{ij} \\
s\endbmatrix_{d_{i}} F^r_i F_j F^s_i = 0,  \qquad if~~ i\neq j ,
 \endalign
$$ where for integers $N, M,
d\geq 0$, we define
$$
[N]_d! = \prod\limits^N_{a=1} \frac{q^{da}-q^{-da}}{q^{d}-q^{-d}},
\qquad \bmatrix M+N \\ N\endbmatrix_d =
\frac{[M+N]_d{!}}{[M]_d{!}[N]_d{!}}.
\tag 1.2 $$

The quantum group
$U_q(\hat {\frak g})$ has a noncocommutative Hopf algebra structure
 with comultiplication $\Delta$,
antipode $S$ and counit $\varepsilon$ defined by
$$\align
\Delta (E_i) &= E_i\otimes 1 + K_i^{}\otimes E_i ,\\
\Delta (F_i) &= F_i\otimes K_i^{-1} + 1\otimes F_i ,\\
\Delta (K_i )&= K_i\otimes K_i,  \tag 1.3\\
S(E_i) &= -K_i^{-1}E_i,\qquad S(F_i)=-F_iK_i^{},\qquad S(K_i)=K_i^{-1} ,\\
\varepsilon (E_i) &= \varepsilon (F_i) = 0,  \qquad \varepsilon (K_i) =1.
 \endalign
$$
\endproclaim

We define an automorphism $D_z$ of $U_q(\hat {\frak g})$ as
 $$D_z(E_0)= z E_0, D_z(F_0)=z^{-1}F_0,  $$
 and $D_z$ fixes   all other
generators.  We also define the map $\Delta_z(a)= (D_z \otimes
id)\Delta(a)$ and $\Delta'_z(a)= (D_z \otimes id)\Delta'(a)$, where
$a\in U_q(\hat {\frak g})$ and $\Delta'$ denotes the opposite
comultiplication.

Let $d$ be an operator such that $d$ commutes with all other elements
but have the relation that $$[d, E_0]=E_0, [d, F_0]=-F_0.   $$
It is clear that action of $D_z$ is equivalent to conjugation by
$z^d$.

For the algebra generated by ${U_q(\hat{\frak g})}$ and
the operator $d$, which we denote by
$U_q(\tilde {\frak g})$, from the theory of Drinfeld\cite{D2},
we know that it has a universal R-matrix $\bar {\frak R}$.

\proclaim {Proposition 1.1}\cite{D2} There exists an element $\bar
{\frak R}$ in $U_q^+(\tilde {\frak g})\hat \otimes  U_q^-(\tilde {\frak g})$
 such that $\bar {\frak R}$ satisfies the properties:
$$\gather
\bar {\frak R}\Delta(a)=\Delta^{\text{op}}(a)\bar {\frak R},\\
(\Delta\otimes \text{id})(\bar {\frak R})=\bar {\frak R}_{13}\bar
 {\frak R}_{23},
\tag 1.4 \\
(\text{id}\otimes\Delta)(\bar {\frak R})=
\bar {\frak R}_{13}\bar {\frak R}_{12},
  \endgather
$$
where $a\in U_q({\frak g})$,   $\Delta^{op}$ denotes the opposite
comultiplication, $\bar
{\frak  R}_{12}=\sum_ia_i\otimes b_i\otimes 1=\bar {\frak R}\otimes 1,
\bar {\frak R}_{13}=\sum_ia_i\otimes 1\otimes b_i,
\bar {\frak R}_{23}=\sum_i1\otimes a_i\otimes b_i=
1\otimes \bar {\frak R}. $
\endproclaim

Here $U_q^+(\tilde {\frak g})$ is the subalgebra generated by
$E_i$, $K_i$ and $d$ and $U_q^-(\tilde {\frak g})$ is the subalgebra
generated by $F_i$, $K_i$ and $d$.

Let $$\bar{\frak R}(z)= (D_z\otimes {\text id})\bar {\frak R}. $$

Let $C$ be  the central element corresponding to the central
extension of the quantum affine algebra.
 Let  ${\frak R}(z)$=$ q^{-d\otimes C -C\otimes d}\bar
{\frak R}(z)$, then
we have
\proclaim {Corollary  1.1} ${\frak R}(z)$
$\in$ $U_q(\hat{\frak  b^+})\hat  \otimes
 U_q(\hat {\frak  b^-})\otimes
{\Bbb C} [[z]]$, such that
$$ {\frak R}(z) \Delta_z (a)= (D^{-1}_{q^{C_2}} \otimes
D^{-1}_{q^{C_1}}) \Delta'_z (a) {\frak R}(z),$$
$$ (\Delta \otimes I){\frak R}(z)= {\frak R}_{13}(zq^{C_2}){\frak R}_{23}(z),$$
$$(I\otimes\Delta){\frak R}(z)= {\frak R}_{13}(zq^{-C_2}){\frak R}_{12}(z),
\tag 1.5 $$
$${\frak R}_{12}(z){\frak R}_{13}(zq^{C_2}/w){\frak R}_{23}(w)=
{\frak R}_{23}(w){\frak R}_{13}(zq^{-C_2}/w){\frak R}_{12}(z).  $$

Here $C_1=C\otimes 1$, $C_2=1\otimes C$,  $U_q(\hat{\frak  b^+})$
is the subalgebra generted by $E_i, K_i$ and  $U_q( {\frak  b^-})$
is the subalgebra generrated by $F_i,K_i$.
\endproclaim

Let
$${\frak R} = (D_z^{-1} \otimes 1) {\frak R}(z). $$

\proclaim {Corollary  1.2}
Let ${\frak C}= ((D^{-1}_{q^{C_2}} \otimes
D^{-1}_{q^{C_1}}) {\frak R}_{21})
{\frak R}$. Then
$${\frak C} \Delta (a)=
((D^{-2}_{q^{C_2}} \otimes
D^{-2}_{q^{C_1}}) \Delta (a)){\frak C}. \tag 1.6  $$
\endproclaim

{\bf Proof.} From the property of ${\frak R}(z)$, we know that
$$ {\frak R} \Delta (a)= (D^{-1}_{q^{C_2}} \otimes
D^{-1}_{q^{C_1}}) \Delta' (a) {\frak R}\tag 1.7   $$
Thus we have
$${\frak R}_{21} \Delta' (a)= ((D^{-1}_{q^{C_2}} \otimes
D^{-1}_{q^{C_1}}) \Delta (a)) {\frak R}_{21}, \tag 1.8   $$
and
$$((D^{-1}_{q^{C_2}} \otimes
D^{-1}_{q^{C_1}}))({\frak R}_{21} \Delta' (a))=
((D^{-1}_{q^{C_2}} \otimes
D^{-1}_{q^{C_1}}))
( ((D^{-1}_{q^{C_2}} \otimes
D^{-1}_{q^{C_1}}) \Delta (a)) {\frak R}_{21}). \tag 1.9   $$
Thus we get that
$$(D^{-1}_{q^{C_2}} \otimes
D^{-1}_{q^{C_1}} {\frak R}_{21}) (D^{-1}_{q^{C_2}} \otimes
D^{-1}_{q^{C_1}}\Delta' (a))= (D^{-2}_{q^{C_2}} \otimes
D^{-2}_{q^{C_1}} \Delta (a))(D^{-1}_{q^{C_2}} \otimes
D^{-1}_{q^{C_1}} {\frak R}_{21}). \tag 1.10   $$
Thus we obtain the proof.

We note that $(q^{C_2}\otimes q^{C_1})$
is invariant under the permutation.

This proposition  shows that the action of
 ${\frak C}=((D^{-1}_{q^{C_2}} \otimes
D^{-1}_{q^{C_1}}) {\frak R}_{21})
{\frak R}$ on a tensor product of two modules
 is  an  intertwiner which, however,   shifts the  actions of
$E_0$ and $F_0$ by  the constants
$q^{\mp 2C_2}\otimes q^{\mp 2C_1}$ respectively.

We should also notice that $D_z\otimes D_z$ acts invariantly on
${\frak R}$. Let ${\frak R}_{21}(z)=(D_z\otimes 1){\frak R}_{21}.$
Note that ${\frak R}_{21}(z)$ is not equal to $P({\frak R}(z))$, where
$P$ is the permutation operator.

Let $V$ be a finite dimensional   representation of
$U_q(\hat{\frak g})$.

Let $${\bar L}^{+}(z)= ({\text id} \otimes
\pi_{V})( {\frak R}_{21}(z)),\tag 1.11a  $$
$${\bar L}^{-^{-1}}(z)= ({\text id} \otimes
\pi_{V}) {\frak R}_{}(z^{-1}), $$
$${\frak L}^{+}(z)= (\pi_{V} \otimes{\text id})
( {\frak R}^{-1}_{}(z)), $$
$${\frak L}^{-^{-1}}(z)= (\pi_V \otimes
{\text id}) {\frak R}^{-1}_{21}(z^{-1}). $$

We have that $${\bar L}^{+}(z) P {\frak L}^{+}(z^{-1})=1, \ \
{\bar L}^{-^{-1}}(z) P {\frak L}^{-^{-1}}(z^{-1})=1, \tag 1.11b $$ where
$P$ is the permutation operator.

$U_q(\hat{\frak g})$ as an algebra
 is  generated  by operator entries of
 ${\bar L}^{+}(z)$ and
${\bar L}^{-^{-1}}(z)$,  and it is also  generated by operator entries
of  ${\frak L}^{+}(z)$ and
${\frak L}^{-^{-1}}(z)$.
${\bar L}^\pm(z)$ are   used in \cite{FR} to obtained
q-KZ equation.

 Let
$${\bar L}({z})= ({\text id} \otimes \pi_{V})((1\otimes
D^{-1}_{q^{C}})
 {\frak R}_{21}({z})){\frak R}(z^{-1}),\tag 1.12a  $$
$$({D^{}_z \otimes 1}) {\bar L}({z})=\bar  L = ({\text id} \otimes
\pi_{V})((1\otimes D^{-1}_{q^{C}} ) {\frak R}_{21}){\frak R},  $$
and $${\frak L}({z})= (\pi_{V} \otimes{\text id})
(( D^{}_{q^{C}} \otimes 1 ) {\frak R}^{-1}
({z})){\frak R}^{-1}_{21}(z^{-1}),\tag 1.12b  $$
$$(1\otimes D^{}_z) {\frak L}({z})= {\frak L} =(\pi_{V}
\otimes{\text id})
((1\otimes D^{}_{q^{C}} ) {\frak R}^{-1}){\frak R^{-1}_{21}}.  $$
\proclaim{Proposition 1.2}
$$\eqalign{
\bar R(\frac{z}{w}) {\bar L}^{\pm}_1(z) {\bar L}^{\pm}_2(w) &= {\bar
L}^{\pm}_2(w) {\bar L}^{\pm}_1(z)
\bar R(\frac{z}{w}),\cr
\bar R(\frac{zq^{-C}}{w}) {\bar L}^+_1(z) {\bar L}^-_2(w) &= {\bar L}^-_2(w)
{\bar L}^+_1(z)
\bar R(\frac{zq^{C}}{w})
,}
$$
$$\bar R(z/w)\bar L_1(z)\bar R(zq^{2C}/w)^{-1}\bar L_2(w)=
\bar L_2(w)\bar R(zq^{-2C}/w)\bar L_1(z)\bar R(z/w)^{-1},\tag 1.13a$$
$${\bar L}({z})(id \otimes \pi_V)\Delta(a)= ((1\otimes D^{-2}_{q^C}) \Delta(a))
{\bar L}(z), $$
and
$$\eqalign{
\bar R(\frac{z}{w}) {\frak L}^{\pm}_1(z)^{-1}
 {\frak L}^{\pm}_2(w)^{-1} &= {\frak L}^{\pm}_2(w)^{-1}
{\frak L}^{\pm}_1(z)^{-1}
\bar R(\frac{z}{w}),\cr
\bar R(\frac{zq^{-C}}{w}) {\frak L}^+_1(z)^{-1} {\frak L}^-_2(w)^{-1} &=
 {\frak L}^-_2(w)^{-1} {\frak L}^+_1(z)^{-1}
\bar R(\frac{zq^{C}}{w})
,}  \tag 1.13b
$$
$${\frak L}_1(z)\bar R(zq^{-2C}/w)^{}{\frak L}_2(w)
\bar R(z/w)^{-1}=  \bar R(z/w)
{\frak L}_2(w)\bar R(zq^{2C}/w)^{-1}{\frak L}_1(z),$$
$${\frak L}({z})(id \otimes \pi_V)\Delta(a)= (D^{}_{q^{2C}}\otimes 1
\Delta(a)){\frak L}(z), $$
Here $\bar R(z/w)$ is the image  of  ${\frak R}(z/w)$ on $V\otimes
V$.
\endproclaim

We name ${\bar L}(z)$ and  ${\frak L}(z)$ universal Casimir operator
and inverse universal Casimir operator of the quantum algebra  respectively.

The construction of
a  representation  $U_q({\hat{\frak g}})$  is equivalent to finding
 a   specific realization  of  the
operator  ${\bar L}(z)$ or ${\frak L}(z)$, which plays the same role as $L$
 in the case of $U_q({\frak g})$ in \cite{DF2}.
Naturally we would like to find a way to build this ${\bar L}(z)$
or ${\frak L}(z)$ out of the
intertwiners, just as in the case of spinor and oscillator
representations
of the  quantum groups of types $A$, $B$, $C$ and $D$ in \cite{DF2}.

Let $V_{\lambda,k}$ and $V_{\lambda_1,k}$
be two  highest weight representations  of
$U_q(\hat{\frak g})$ with highest weight $\lambda$ and
$\lambda_1$ and the center $C$
acting  as a  multiplication by  a number $k$.

 Let $\Psi$ be an  intertwiner:
$$\Psi:V_{\lambda_1,k} \longrightarrow V_{\lambda,k} \otimes V,$$
$$\Psi (x)= \Psi_1(x)\otimes e_1+...+ \Psi_n(x)\otimes e_n,\tag 1.14  $$
where $x\in V_{\lambda_1, k}$ and $e_i$ is a  basis for $V$.
Let $\Psi^*$ be an  intertwiner:
$$\Psi^*:V_{\lambda,k} \longrightarrow V_{\lambda_1,k} \otimes V^*,$$
$$\Psi^* (x)= \Psi_1^*(x)\otimes e_1^*+...+ \Psi_n^*(x)\otimes
e_n^*,\tag 1.15  $$
where $V^*$ is the left dual representation of $V$ of $\hat U_q{({\frak
g}})$, $x\in V_{\lambda,k}$ and $e_i^*$ is a   basis for $V^*$.
Let us identify $V\otimes V^*$ with End$(V)$.

$(\Psi \otimes 1)\Psi^*=
\Sigma \Psi_i \Psi^*_j\otimes e_i\otimes
e^*_j$ gives a map
 $$(\Psi \otimes 1)\Psi^*:V_{\lambda,k} \longrightarrow V_{\lambda,k} \otimes
V\otimes
V^*.$$

Let $\tilde L\in$
End$(V_{\lambda,k}) \otimes$ End$(V)$:
$$\tilde L= (\tilde L_{ij})= ( (D_{q^{-2k}}\Psi_i) \Psi_j^*), \tag
1.16 $$ where $(D_{q^{-2k}}\Psi_i)$ means   shifting  the evaluation
representation by the  constant $q^{-2k}$. Here we need the
assumption that  $\tilde L$ is
well defined, which we define as follows:

On $V_{\lambda,k}$ and $V_{\lambda_1,k}$,
 we define a grading by defining the action of
$E_0$ lifting the degree of an element by 1 and
$F_0$ lowering the degree of an element by 1.
Then $\Psi_i$ can be written as $ \Sigma_{\Bbb Z} \Psi_i(n)$,
where $\Psi_i(n)$ shifts the degree of an element by $n$.
Let $\Psi_i^{a,b}=\Sigma_{a}^b \Psi_i(n)$, $a<b$.
  By $\tilde L$ is well defined, we mean that any homogeneous
component of  the image of  an element in $V_{\lambda,k}$ under the
composite action of  $\Psi_j^*$ and $\Psi_i^{a,b}$
converges when $a$ goes to negative infinity and
$b$ goes to positive infinity.

\proclaim{Proposition 1.3}$$
\tilde L \Delta(a)=((1\otimes D_{q^{-2k}})\Delta(a))\tilde L, \tag
1.17 $$
where $a$ is an element in $U_q(\hat{\frak g})$.
\endproclaim

 Let $\Phi$ be an intertwiner:
$$\Phi:V_{\lambda_1,k} \longrightarrow V\otimes V_{\lambda,k} ,$$
$$\Phi (x)=e_1\otimes  \Phi_1(x)+...+ e_n\otimes \Phi_n(x),\tag 1.18  $$
where $x\in V_{\lambda_1, k}$ and $\{ e_i\}$ is the  basis for $V$.
Let $\Phi^*$ be an  intertwiner:
$$\Phi^*:V_{\lambda,k} \longrightarrow ^*V\otimes  V_{\lambda_1,k} ,$$
$$\Phi^* (x)= e_1^*\otimes\Phi_1^*(x)+...+e_n^*\otimes  \Phi_n^*(x)
,\tag 1.19  $$
where $^*V$ is the right dual representation of $V$ of $\hat U_q{({\frak
g}})$, $x\in V_{\lambda,k}$ and $e_i^*$ is the  basis for $^*V$.

 By the right dual representation
of $V$ of $\hat U_q{({\frak
g}})$, we mean the action of $\hat{U_q({\frak
g}})$ on the dual space given by $<av',v>=<v',S^{-1}(a)v>$, for
$a\in \hat U_q{({\frak g}})$, $v\in V$ and $v' \in$ $ ^*V$.

$(1\otimes \Phi)\Phi^*=
\Sigma \Phi_i \Phi^*_j\otimes e^*_j\otimes
e_i$ gives a map
 $$\Phi:V_{\lambda,k} \longrightarrow  ^*V \otimes V \otimes
V_{\lambda,k} .$$

Let us identify $^*V\otimes V$ with End$(V)$
by the  map  which maps  the first two components of
 $V\otimes ^*V\otimes V$ to $\Bbb C$ and fix the last component.

Let $\tilde {\frak L}\in$
 End$(V)\otimes$  End$(V_{\lambda,k})$:
$$\tilde {\frak L}=
 (\tilde {\frak L}_{ij})= ((D_{q^{2k}}\Phi_j)\Phi_i^*  ), \tag
1.20 $$ where $(D_{q^{2k}}\Phi _i)$ means shifting  the evaluation
representation by constant $q^{2k}$ and we assume $\tilde {\frak L}$ is
well defined.
\proclaim{Proposition 1.5}
 $$\tilde {\frak L} \Delta (a)=  (( D_{q^{2k}}\otimes 1) \Delta(a)) \tilde
{\frak L}\tag 1.21    $$
\endproclaim

The key idea  is to identify   $\tilde L$ with
 $L$ or $\tilde {\frak L}$ with ${\frak L}$
 to obtain representations out of intertwiners.

The intertwiners for the affine quantum groups are extensively
studied by Kyoto school. They  connected the XXZ model in statistical
mechanics with the structures of the representation of
quantum affine algebras via the intertwiners.
Meanwhile Jimbo, Miki, Miwa and Nakayashiki \cite{JMMN} and Koyama \cite{K}
 worked on the bosonlization of the intertwiners of
quantum affine algebras. These results, in some sense,  is trying to
 obtain the fermions out of  bosons.

The results  in \cite {DF} enable us to  obtain
the explicit  quantum boson-fermion  correspondence via
 Gauss decomposition of $L^\pm(z)$.
Our construction
 also brings a  different but conceptional understanding
to the  classical boson-fermion correspondence.

 On the other hand,    Miki\cite{M},
Foda, Iohara, Jimbo, Kedem and Yan \cite{FIJKMY} presented another
idea to construct realizations  of ${\bar L}^{\pm}$ by composition of
$\Psi$ and $\Phi^*$.

In a subsequent paper,
we will use the idea
 in this paper to define the corresponding quantum Clifford algebras,
quantum Weyl algebras  and  construct spinor and oscillator
representations of $U_q(\hat {\frak o}(2n))$ and $U_q(\hat {\frak
sp}(2n))$.
 For the spinor representation of
$U_q(\hat {\frak o}(2n))$, only the
bosonic  realization \cite{B} is available.
 We expect to  build the boson-fermion
correspondence in a similar  but simpler  way.

\documentstyle{amsppt}
\subheading{ 2. Realizations of  Affine Quantum Groups}

The Drinfeld-Jimbo definition of quantum groups  by generators and
relations is valid for an arbitrary generalized Cartan matrix.
In particular, the choice of extended Cartan matrix of type  $A_n^
{(1)}$ yields the quantum affine algebra $U_q(\widehat{{\frak{sl}}(n)})$.
Drinfeld found in \cite{D4} another realization of the
quantum affine algebras, which to a certain degree plays the
role of loop algebra realization in the undeformed case.
 We extend Drindeld's construction to the quantum affine
 algebra $U_q(\widehat{{\frak{gl}
}(n)})$.

\proclaim{Definition 2.1}  $U_q(\widehat{{{\frak{gl}}}(n)})$ is an
associative algebra with unit 1 and generators\linebreak
$$\{X^{\pm}_{ik},k^+_{jl}, k^-_{jm}, q^{\pm \frac12c}\mid
i=1,\dots, n-1,j=1,\dots,n, k\in\bold Z, l\in\bold Z_+, m\in -\bold
Z_+\} \tag 2.1 $$
 satisfying relations in terms of the following
generating functions in a formal variable $z$: $$\eqalign
{X^{\pm}_i(z)  = \sum_{k\in\bold Z} X^{\pm}_{ik}z^{-k} \cr k^+_j(z)  =
\sum_{l\in
\bold Z_{+}} k^+_{jl} z^{-l} \cr k^-_j(z)  =
\sum_{m\in \bold -Z_{+}} k^-_{jm}z^{-m}.  }\tag 2.2 $$
The generators $q^{\pm \frac 12c}$ are central and mutually
inverse. The other relations are:
 $$ k^+_{j0} k^-_{j0} = k^-_{j0} k^+_{j0} = 1
.  \tag 2.3 $$
$$\eqalign{k^{\pm}_i(z) k^{\pm}_j(w) & = k^{\pm}_j(w) k^{\pm}_i(z)  \cr
k^+_i(z) k^-_i(w) & = k^-_i(w) k^+_i(z)
.}$$
$$\eqalign{\frac{z_{\mp}-w_{\pm}}{z_{\mp}q^{-1}-w_{\pm}q^{}} k^{\mp}_i(z)
k^{\pm}_j(w) & = k^{\pm}_j(w) k^{\mp}_i(z)
\frac{z_{\pm}-w_{\mp}}{z_{\pm}q^{-1}-w_{\mp}q^{}}, \quad\text{if}\quad
j>i
.}$$
$$\eqalign{\cases k^{\pm}_i(w)^{-1} X^{\pm}_j(z) k^{\pm}_i(z) & =
X^{\pm}_j(z)\qquad\text{if}\quad i-j\leq -1 \quad\text{or}\\
k^{\pm}_i(w)^{-1} X^{\mp}_j(z) k^{\pm}_i(w)&  = X^{\mp}_j(z)
\qquad\quad i-j\geq 2\endcases
. }$$
$$\eqalign{k^{\pm}_i(z)^{-1}X^-_i(w) k^{\pm}_i(z)  &=
\frac{z_{\mp}q^{-1}-wq^{}}{z_{\mp}-w} X^-_i(w) \cr
k^{\pm}_{i+1}(z)^{-1}X^-_i(w) k^{\pm}_{i+1}(z)&  =
\frac{z_{\mp}q^{}-wq^{-1}}{z_{\mp}-w} X^-_i(w)  \cr
k^{\pm}_i(z) X^+_i(w) k^{\pm}_i(z)^{-1} & =
\frac{z_{\pm}q^{-1}-wq^{}}{z_{\pm}-w} X^+_i(w) \cr
k^{\pm}_{i+1}(z) X^+_i(w) k^{\pm}_{i+1}(z)^{-1} & =
\frac{z_{\pm}q^{}-wq^{-1}}{z_{\pm}-w}X^+_i(w)
.}$$
$$\eqalign{ (zq^{}-wq^{-1})X^-_i(z) X^-_i(w) & = X^-_i(w) X^-_i(z)
(zq^{-1}
-wq^{})
  \cr (zq^{-1}-wq^{}) X^+_i(z) X^+_i(w) & = X^+_i(w) X^+_i(z) (zq^{}-wq^{-1})
.}\tag 2.4 $$
$$\eqalign{ (z-w)X^+_i(z) X^+_{i+1}(w)&  =
(zq^{-1}-wq^{}) X^+_{i+1}(w)X^+_i(z)
\cr (zq^{-1}-wq^{})X^-_i(z) X^-_{i+1}(w)
& = (z-w) X^-_{i+1}(w) X^-_i(z)\cr
 [X^{\pm}_i(z), X^{\pm}_j(w)] & = 0,
\qquad\text{for}\quad A_{ij}=0}. $$
$$ [X^+_i(z), X^-_j(w)]  = $$
$$ (q-q^{-1})\delta_{ij}\{\delta (zw^{-1}
q^{-c})k^-_{i+1}(w_-)k^-_i(w_-)^{-1}  - \delta(zw^{-1}q^c)
k^+_{i+1}(z_-)k^+_i(z_-)^{-1}\}.
$$
$$\eqalign{
\{X^{\pm}_i(z_1)X^{\pm}_i(z_2) X^{\pm}_j(w)  - (q+q^{-1})
X^{\pm}_i(z_1) X^{\pm}_j(w) X^{\pm}_i(z_2) +
\cr
X^{\pm}_j(w) X^{\pm}_1(z_1) X^{\pm}_i(z_2)\}  + \{z_1\leftrightarrow
z_2\} = 0, \cr
\text{for}\quad A_{ij}=-1
.}$$
Here $z_{\pm} =  q^{\mp \frac{c}{2}}z$,
$A_{ij}$ are the entries of the Cartan matrix for
${{\frak{sl}}(n)}$, and
$$\delta(z) = \sum_{n\in\bold Z}z^n.  \tag 2.5 $$

Drinfeld's  realization of the subalgebra
$U_q(\widehat{{\frak{sl}}(n)})$ is
given by
$x^{\pm}_i(z)=$ $(q-q^{-1})^{-1}
X^{\pm}_i(zq^i)$, $\psi_i(z) =$$k^-_{i+1}(zq^i) k^-_i(zq^i)^{-1}$ and
$\varphi_i(z)=k^+_{i}(zq^i)k^+_{i+1}(zq^i)^{-1}$.
 \cite {DF}.

\endproclaim

 Drinfeld \cite{D3}  stated  that
the  algebra  $U_q(\widehat{{\frak{sl}}(n)})$  is isomorphic to
the one constructed by generators and relations with
 extended  Cartan matrix of
 type ${\hat{A}}_{n-1}$.

The third and the fourth formulas above are defined on the
completion of the tensor algebra.

Let  $R(z)$ be an element of End$(\Bbb C^n\otimes \Bbb C^n)$ defined by
$$\eqalign{
R(z) &= \sum^n_{i=1}E_{ii}\otimes E_{ii} + \sum^n\Sb i\neq
j\\ i,j=1\endSb E_{ii}\otimes E_{jj}\frac{z-1}{q^{-1}z
-q}\cr
&+ \sum^n\Sb i<j\\ i,j=1\endSb E_{ij}\otimes
E_{ji}\frac{z(q^{-1}-q)}{z q^{-1}-q} + \sum^n\Sb i>j\\
i,j=1\endSb E_{ij}\otimes E_{ji} \frac{(q^{-1}-q)}{q^{-1}z-q}
} \tag 2.6
$$
where $q,z$ are formal variables. Then
$R(z)$ satisfies the  Yang-Baxter equation with a parameter:
$$
R_{12}(z) R_{13}(z/w) R_{23}(w) = R_{23}(w)
R_{13}(z/w) R_{12}(z),
\tag 2.7 $$
 and $R$ is unitary, namely
$$
R_{21}(z)^{-1} = R({z^{-1}}).
\tag 2.8
$$

The definition of $R(z)$ implies
$$\eqalign{
\lim\limits_{z\to 0}R(z) &= \sum^n_{i=1} E_{ii}\otimes
E_{ii} +  \sum^n\Sb i\neq j\\ i,j=1\endSb
qE_{ii}\otimes E_{jj} +q^{-1} (q^{}-q^{-1})\sum^n\Sb i>j\\ i,j=1\endSb
E_{ij}\otimes E_{ji}\cr
&= q^{-1}R_{21}^{},\cr
\lim\limits_{z\to \infty}R(z) &= \sum^n_{i=1} E_{ii}\otimes
E_{ii} +  \sum^n\Sb i\neq j\\ i,j=1\endSb
q^{-1}E_{ii}\otimes E_{jj} +q^{} (q^{-1}-q^{})\sum^n\Sb i<j\\ i,j=1\endSb
E_{ij}\otimes E_{ji}\cr
&=q^{} R^{-1} ,
}\tag 2.9
$$
where $R$ is the R-matrix of
$U_q({\frak gl}(n))$  on $\Bbb C^n \otimes \Bbb C^n$\cite{J3}.

Faddeev, Reshetikhin and Takhtajan defined a Hopf  algebra using this
element $R(z)$.  Reshetikhin and
Semenov-Tian-Shansky obtained a central extension of an  algebra,
which is defined in the same way but with a R-matrix of a scalar
multiple of $R(z)$ above. However these two algebras are isomorphic,
which is implied in the proof of the main theorem in \cite{DF}.

We will denote by $U(\tilde R)$ the algebra with central extension
defined by the $R(z)$ above.
The central extension is incorporated in shifts of the parameter $z$
in $R(z)$.

\proclaim{Definition 2.2}  $U(\tilde R)$ is an associative algebra
with generators $\{l^{\pm}_{ij}[\mp m], m\in{\bold Z_+ \setminus  0}
\quad\text{and}\mathbreak
l^+_{ij}[0], l^-_{ji}[0], 1\leq j\leq i\leq n\}$.  Let
$l^{\pm}_{ij}(z) = \sum\limits^{\infty}_{m=0} l^{\pm}_{ij}[\pm
m] z^{\pm m}$, where $l^+_{ij}[0] = l^-_{ji}[0] = 0$, for $1\leq
i<\leq n$.  Let $L^{\pm}(z) = (l^{\pm}_{ij}(z))^n_{i,j=1}$.
Then  the defining  relations are the following:
 $$
l^+_{ii}[0] l^-_{ii}[0] = l_{ii}[0] l^+_{ii}[0]=1,
$$
$$\eqalign{
R(\frac{z}{w}) L^{\pm}_1(z) L^{\pm}_2(w) &= L^{\pm}_2(w) L^{\pm}_1(z)
R(\frac{z}{w}),\cr
R(\frac{z_+}{w_-}) L^+_1(z) L^-_2(w) &= L^-_2(w) L^+_1(z)
R(\frac{z_-}{w_+})
,}\tag 2.10
$$
where $z_{\pm} = zq^{\mp  \frac{c}{2}} $.
For the first formula of (2.10), we can expand
$R(\frac{z}{w})$   in formal power series of either
 $\frac{z}{w}$ or $\frac{w}{z}$, but for
the second  formula of (2.10), it is
 expanded in the formal power series of  $\frac{z}{w}$.

$U(\tilde R)$ is a Hopf algebra:  its coproduct is defined by
$$\eqalign{
\Delta' L^{\pm}(z) = L^{\pm}(zq^{\mp(1\otimes\frac{c}{2})})
&\dot{\otimes} L^{\pm}(zq^{\pm(\frac{c}{2}\otimes 1)})\cr
\text{or}\qquad \Delta'(l^{\pm}_{ij}(z)) = \sum^n_{k=1}
l^{\pm}_{ik}(zq^{\mp(1\otimes\frac{c}{2})}) &\otimes
l^\pm_{kj}(zq^{\pm(\frac{c}{2}\otimes 1)}),
}\tag 2.11
$$
and its antipode is
$$
S(L^{\pm}(z)) = L^{\pm}(z)^{-1}.\tag 2.12
$$
Note that the invertibilty of $L^{\pm}(z)$ follows from the properties  that
  $l^{\pm}_{ii}$ are invertible and $L^{+}(0)$ and $L^{-}(0)$
 are upper triangular and
lower triangular operator-entried matrices  respectively.

\endproclaim
\vskip .25in

\proclaim{ Theorem 2.1}
$L^{\pm}(z)$ have the following unique decompositions:
$$\align
L^{\pm}(z) = &\pmatrix 1& & &  0\\
 e^{\pm}_{2,1}(z) &\ddots & &\\ e^{\pm}_{3,1}(z)
&\ddots &\ddots\\ \vdots &\ddots &\ddots &\ddots\\ e^{\pm}_{n,1}(z)
&\hdots &e^{\pm}_{n,n-1}(z) &e^{\pm}_{n-1,n}(z) &1\endpmatrix
\pmatrix k^{\pm}_1(z) & & &  0 \\ &\ddots\\  & &\ddots\\ &\\0 & &
&k^{\pm}_n(z)\endpmatrix\times \tag 2.13\\  \\
&\pmatrix 1 & f^{\pm}_{1,2}(z) & f^{\pm}_{1,3}(z) &\hdots
&f^{\pm}_{1,n}(z)\\ &\\ &\ddots &\ddots &\ddots &\vdots\\ & & &
&f^{\pm}_{n-1,n}(z)\\ &\\0  & & & &1\endpmatrix
,\endalign
$$
where
$e^{\pm}_{i,j}(z)$, $f^{\pm}_{j,i}(z)$ and
$k^{\pm}_i(z)$($i>j$)
 are elements in $U(\tilde{R} )$  and $k^{\pm}_i(z)$ are invertible.  Let
$$\eqalign {
X^-_i(z) &= f^+_{i,i+1}(z_+) - f^-_{i,i+1}(z_-),\cr
X^+_i(z) &= e^+_{i+1,i}(z_-) - e^-_{i+1,i}(z_+),}\tag 2.14
$$ then
 $q^{\pm\frac12c}, X^{\pm}_i(z), k^{\pm}_j(z), i=1,\dots,n-1,
j=1,\dots,n$ satisfy the relations (2.3), (2.4)
 of $U_q(\widehat{{{\frak{gl}}}(n)})$. The homomorphism
$$ M  :  U_q(\widehat{{{\frak{gl}}}(n)})
\longrightarrow  U(\tilde{R}) \tag 2.15 $$
defined by (2.14) is an isomorphism.
\endproclaim

Since  $k^{\pm}_i(z)$ are invertible,
the elements $e^{\pm}_{i,j}(z)$, $f^{\pm}_{j,i}(z)$ and
$k^{\pm}_i(z)$($i>j$)
 are uniquely expressed in terms of the matrix coefficients of
$L^\pm(z)$. In view of this analogy
we call (2.13)   Gauss
decomposition
of
$L^{\pm}(z)$.

We note that as  a corollary of Theorem 2.1  and Definition 2.1,
one gets a realization of $ U_q(\widehat{{{\frak{sl}}}(n)})$ as a
subalgebra of $U(\tilde R)$.

Let $L(z)$= $L^+(zq^{-c})(L^-(z))^{-1}$, then
$$ R(z/w)L_1(z) R(zq^{2c}/w)^{-1}L_2(w)=
L_2(w) R(zq^{-2c}/w)L_1(z) R(z/w)^{-1}. \tag 2.16 $$

\documentstyle{amsppt}
{\bf
3. Spinor
representation of
$\hat {\frak gl}(n)$ and
 Boson-Fermion Correspondence}

For $\hat {\frak gl}(n)$, we  assume $n>2$. This restriction is
to avoid some  nonessential but tedious
 complication caused by the selfdual property of
the standard representation  of $ {\frak sl}(2)$ on $\Bbb
C^2$, which can  be resolved.

 \proclaim{Definition 3.1}
The Affine  Clifford algebra is an associative
algebra   generated by
$a_i(m)$ and $a^*(m)$, ($i=1,..,n;
m \in \Bbb Z$) with the commutation relations:
$$ a_i(m)a_j(l)+a_j(l)a_i(m)=0,$$
$$a^*_i(m)a^*_j(l)+a^*_j(l)a^*_i(m)=0, $$
$$a_i(m)a^*_j(l)+a^*_j(l)a_i(m)=\delta_{ij}\delta_{m,-l}1.\tag 3.1  $$
\endproclaim

\proclaim{Definition 3.2}The affine Heisenberg algebra is
an associative algebra generated by
$h_i(m)$, $m\neq 0, i=1,...,n$ and $m\in \Bbb Z$. The relations are
$$ h_i(m)h_j(l)-h_j(l)h_i(m)=0,l\neq -m, $$
$$ h_i(m)h_j(-m)-h_j(-m)h_i(m)=m\delta_{ij},m\neq 0. \tag 3.2 $$
\endproclaim
We now introduce the  notation  of formal power series. (For detail see
\cite{FLM}.) For a vector space $W$, we denote that
$$ W[[z,z^{-1}]]=\{\Sigma_{m\in \Bbb Z} v_mz^m\|v_m \in W \}, \tag 3.3
$$
 where
 $z$ is  a formal variable.  For example, let $W$ be also
an algebra on $\Bbb C$, $ \delta(z)= \Sigma_{m\in \Bbb Z}z^m$ as defined in
(2.5) will be used frequently.
Given two former power series $g_1(z)=\Sigma g_i(m)z^{-m}$ and
$g_2(z)=g_2(m)z^{-m}$ over  an algebra
$W$,  by the product of
these two formal power series  operators is  well
defined, we mean that the limit of each homogeneous  component of
$\Sigma_{m_1<m<m_2} g_1(m)z^{-m}\Sigma_{l_1<l<l_2} g_2(l)z^{-l}$,
when $l_1,m_1$ go to negative infinity and $m_2, l_2$ go to
positive infinity,  exists. We also define the product as a formal
power series with the corresponding limit as the coefficient
of the corresponding homogeneous component. This complies with the
definition of the composition of intertwiners in Section.
For this paper, we define that the limit of the series
of operators  as the operator whose matrix coefficients are the
limits of the corresponding series of
matrix coefficients  of the series of operators.

For example, let
$f(z)$ be a polymonial of $z$ and $z^{-1}$, then
$f(z)\delta(z)$ is well defined and
 $$f(z)\delta(z)=f(1)\delta(z). \tag 3.4$$

Let $E_{ij}$ be the standard basis of Lie algebra
 ${\frak gl}(n)$. As showed in
\cite{G},  $\hat{\frak gl}(n)$ can be realized as ${\frak gl}(n)\otimes
[x, x^{-1}]\oplus C$, where C is the central element. Let
$E_{ij}(z)= \Sigma  E_{ij}\otimes x^n z^{-n}$ be
 the standard basis of the generating functions of
$\hat{\frak gl}(n)$, where  $z$ is a formal variable.

We define $a_i(z)$, $ a^*(z)$
formal  power series with coefficients in the corresponding
 algebra as
$$ a_i(z)=\Sigma a_i(m)z^{-m},$$
$$ a^*_i(z)=\Sigma a^*_i(m)z^{-m}. \tag 3.5 $$

Let $P$ be an n dimensional lattice  $\Sigma \oplus \Bbb Z
h(0)_i$, with the form that $(h(0)_i, h(0)_j)=\delta_{ij}$.
Let $\Bbb C[ \bar P]$ be
the central extension  of $\bar P$, the group algebra of
$P$, such that
$e^{h(0)_i} e^{h(0)_j}$= $(-1)^{(h(0)_i, h(0)_j)}e^{h(0)_j}e^{h(0)_j}$.
Let $\bar H$ be an  associative  algebra  generated by
$h_i(m)$, $m\ne 0$ and $\Bbb C[ \bar P]$, where $h_i(m)$ and $\Bbb C[ \bar P]$
commute with each other.
Let $h_i(z)=\Sigma h_i(m)z^{-m} + \partial_{h(0)_i}$, where
$\partial_{h(0)_i}$ is the partial differential of $h(0)_i$.

\proclaim{Definition 3.3} Spinor Fock space is the space of the subalgebra
of affine Clifford algebra
generated  by $a_i(-m)$, $a^*(-m)$, $m>0$ and $a^*_i(0)$.
 \endproclaim

\proclaim{Definition 3.4} Oscillator Fock space
is the space of the subalgebra of $\bar H$  generated by $h_i(-n), n>0$ and
$\Bbb
C[\bar P]$.
\endproclaim

In this paper, we define normal ordering  : : as in \cite{F1}.

\proclaim{Proposition 3.1} \cite{F} \cite{FF} \cite{F1}
Let $a_i(z), a^*(z)$ be as defined above and
$E_{ij}(z)$ be the standard basis of $\hat{\frak gl}(n)$.
There is a representations of $\hat{\frak
gl}(n)$ given by the following on the spinor Fock space:
$$E_{ij}(z)\longrightarrow :a_i(z)a^*_j(z) :.\tag 3.6 $$
where : : denotes the  normal ordering.
\endproclaim

\proclaim{Proposition 3.2}\cite{FK} \cite{F1}
Let $h_i(z)$ as defined above,
there is  the  Frenkel-Kac construction of  a representation of
$\hat{\frak gl}(n)$ on the oscillator Fock space  by vertex operators given by
$$E_{ii}(z)\longrightarrow h_i(z),\tag 3.7  $$
$$E_{ij}(z)\longrightarrow :{\text  exp}(\bar h_i(z)-\bar h_j(z)):,
i\neq j, $$
where $\bar h_i(z)$=$\Sigma_{n\neq 0} 1/n h_i(-n) z^n+
\partial_{h(0)_i} { ln}z+h(0)_i$.
This representation  is isomorphic
to the representation of Proposition 3.1  above.
The isomorphism is  given  as:
$$ a_i(z) =:{\text  exp}(\bar h_i(z)):, $$
$$a^*_i(z) =:{\text  exp}(-\bar   h_i(z)):. \tag 3.8$$
Here $: :$ is defined in the standard way as in \cite{F1}.
\endproclaim

The isomorphism can be easily proved. First, we shift the degree
of $ a_i(z)$ and $a^*_i(z)$ by $\pm 1/2$ respectively. Then
 by comparing the character, we get the proof\cite{F1}.
The propositions above give us  the classical
boson-fermion correspondence.

 Let's denote the representation
by $V_{bf}$.
\proclaim{Proposition 3.3}
$V_{bf}$ contains a unique copy  of  the highest weight
representations of $\hat {\frak gl}(n)$ with C acting as a
multiplication of 1,
 the highest weight as any   of the fundamental weights corresponding to
any i-th fundamental weight of  the classical part ${\frak sl}(n)$
of $\hat{\frak sl}(n)$ or  zero weight and the action of
$\Sigma_{i=1,...,n} E_{ii}$ as any  integer $l$ such that
$l=i({\text mod}(n))$.
\endproclaim

 We know that  $\hat{\frak gl}(n)$ $=$
$\hat{\frak sl}(n)$ $\oplus$  $\hat{\frak gl}(1)$. For the above case,
$\Sigma h_i(z)= \Sigma_{i=1,...n} :a_i(z)a^*_i(z) :$ gives us
$\hat{\frak gl}(1)$, which commutes with $\hat{\frak sl}(n)$.

Let $e_i$ be  the standard
  basis of $V= \Bbb C^n$ and $e_i^*$ be the dual basis of the
dual module $V^*$.
Let $V(z)$ and $V^*(z)$ be the
evaluation representation of $\hat{\frak gl}(n)$, where $z$ is
a formal variable and the action of
$E_{ij}\otimes t^n$ of $\hat{\frak gl}(n)$ is given by
the action of $E_{ij}$ of ${\frak gl}(n)$ with $z^n$ scalar multiple.

\proclaim{ Proposition 3.4} Let
$$\Phi^c=\Sigma a^*_i(z) \otimes  e_i, $$
$$\Phi^{c*}=\Sigma a_i(z) \otimes e_i^*.\tag 3.9  $$
Then  $\Phi^c$ and $\Phi^{c*}$
are intertwiners:
$$ \Phi^c: V_{bf} \longrightarrow V_{bf}\otimes V_z,$$
$$ \Phi^{c*}: V_{bf} \longrightarrow V_{bf}\otimes V^*_z.\tag 3.10 $$
\endproclaim

This can be checked by direct calculation. This  observation
is the starting point of our  approach to
the spinor and oscillator representations.
As in the finite dimensional cases \cite{DF2},
we  explain  the hidden structure behind all the constructions above.

Now we will derive the above construction without any specific
realizations. In the next section, we will
present a  parallel q-deformation
of this abstract construction.

We now start from an abstract module
$V_{bf}$=$\Sigma_{l\in \Bbb Z} V^l$,
where $V^l$ is the  highest weight representation
with the $l({\text mod}(n))$-th fundamental weight,
 with central extension 1 and
 the action of $\Sigma E_{ii}$ as integer
$l$. By the 0-th fundamental weight, we mean zero weight.
$V^l$ has infinite many copies of
the highest weight representation of $\hat {\frak sl}(n)$ corresponding to the
the $l({\text mod}(n))$-th fundamental weight of
 ${\frak sl}(n)$ and  with central extension 1.
{}From the abstract representation theory of Kac-Moody algebras, we know
the existance of those modules.

As we explain above, the  representation  of
$\hat {\frak gl}(n)$ can be derived on the  the space
$ V_z= V[z, z^{-1}]$ with the action $\Sigma  E_{ii}$ as $1$.

We refer this notation to \cite{FLM}, where it  is
introduced for vertex operator algebras.
By $\bar V_z$, we mean the set of all vectors in the form of
$\Sigma v_if_i(z)$, where $v_i\in V$ form  a basis of
$V$ and $f_i(z)$ are  formal power series
on $z$ and $z^{-1}$.

Let $V^*$ be the left dual module of $V$. Let $F$ be  an  invariant
vector in $V\otimes V^*$ of ${\frak
gl}(n)$, which is unique up to a scalar multiple.
We normalize it,  such that $F$ is equal to the  identity if it is
 identified
as an element in End$(V)$.

\proclaim { Proposition 3.5 }
$F(z_1,z_2)=\{x| x=Ff(z_1)\delta(z_1/z_2)\}$, where
$f(z_1)$ is a polynomial of $z_1$ and $z_1^{-1}$,
 is an  invariant subspace  of
$\bar V_{z_1}\otimes \bar V^*_{z_2}$.
\endproclaim

\proclaim{Proposition 3.6}
Let $\bar \Psi_z^{l}$ and $\bar \Psi_z^{*l}$ be intertwiners as the following:
$$\bar  \Psi_z^{l}:   V^{l+1}
 \longrightarrow V^{l}\otimes V_z,\tag 3.11  $$
$$ \bar \Psi_z^{*l}:   V^{l}
\longrightarrow V^{l+1}\otimes V^*_z,   $$
such that $\bar \Psi_z^{l}= \Sigma \bar \Psi_i^{l}(-m)\otimes v_iz^m$ and
 $\bar \Psi_z^{*l}= \Sigma \bar \Psi_i^{*l}(-m)\otimes v_iz^m$, where
$\bar \Psi_i^{l}(m)$ and  $\bar \Psi_i^{*l}(m)$ are operators  shifting the
degree by $m$.
Then these  operators   are  unique up to a scalar
multiple.
\endproclaim

However, from \cite{FF}, we know that not all the highest weight vectors
are graded 0, but that most of them have negative grading. We will thus shift
the grading of $V^{l}$, for
$l<0$ or $l>n$ by defining that
the grad of the highest weight vector of $V^l$is equal to
 $-n((\Sigma_0^{|m|-1})
+j(m)$, where $l=mn+j$, $n-1>j>0$. After the shifting, we will denote the
new operators by  $ \Psi_z^{l}$ and $\Psi_z^{*l}$.

We know that the correlation of those operators satisfying
KZ equation\cite{TK} \cite{FR},
the solutions for the KZ equation  in this case
can be obtained. On the other hand, because the
 the space of the  solutions  in this case is one dimensional,
we  can normalize them in such a way that
all $<v_{V^{l-1}},\Psi_{z_2}^{l}\Psi_{z_1}^{l+1}
 v_{V^l}>$(i=0,1...,n-1)
are equal, all  $<v_{V^{l+2}},\Psi_{z_2}^{*l+1}\Psi_{z_1}^{*l}
 v_{V^l}>$(i=0,1...,n-1)
are equal and all $<v_{V^l},\Psi_{z_2}^{l}\Psi_{z_1}^{*l}  v_{V^l}>$
(i=0,1,...,n-1)
are equal to $1/(1-z_1/z_2)F$, where $v_{V^l}$ is the highest
weight vector of $V^l$. More precisely, we can choose the
normalization in the way that, their correlation functions are equal
to the corresponding correlation functions of
$\Phi^c$ and $\Phi^{c*}$.

Let
$\Psi_z=\Sigma _i\oplus \Psi_z^{i}$ and
$\Psi_z^*=\Sigma _i\oplus \Psi_z^{*i}$, then $\Psi_z$ and $\Psi_z^*$
are  intertwining operators  from  $V_{bf}$ to $V_{bf}\otimes  V_{z}$
and
$V_{bf}\otimes  V^*_{z}$ respectively.
For  any   vector $v$  in $V$, we can get
${\Psi_{}}(v)(n)$ in
End$(V_{bf})$ and the same of ${\Psi_{}}^*(v^*)(n)$ for $v^*$ in $V^*$.

\proclaim{Proposition 3.7}
$$(\Psi_{z_2}\otimes I)\Psi_{z_1}+ P' (\Psi_{z_1}\otimes I)\Psi_{z_2}=0,$$
$$(\Psi_{z_2}^{*}\otimes I)\Psi_{z_1}^{*}+P'
(\Psi_{z_1}^{*}\otimes I)\Psi_{z_2}^{*}=0, \tag 3.12 $$
$$(\Psi_{z_2}^{}\otimes I)\Psi_{z_1}^{*}+P'
(\Psi_{z_1}^{*}\otimes I)\Psi_{z_2}^{}-\delta(z_1/z_2)F=0,$$
where $P'$ is the operator  which maps $a_i z_1^m\otimes b_jz_2^l$ to
$b_jz_2^l\otimes a_i z_1^m$.
\endproclaim

{\bf Proof.}  We can first prove the relations on the level of
the correlation of the highest weight vector. Because
 the left hand sides   of the formulas  above are intertwiners, we
can prove the equality above for correlation functions
of any two vectors. Thus we finish the proof.

\proclaim{Theorem 3.1}The algebra generated by
${\Psi_{}}(v)(m)$ and ${\Psi_{}}^*(v^*)(m)$ is isomorphic to  the
Clifford algebra of Definition 3.2.
        \endproclaim

The proof is the same as that for Proposition 3.2.

  Let ${\frak g}$ be a simple Lie algebra over $\Bbb C$ of type
$A_n$,  $B_n$, $C_n$  or  $D_n$. Let $e_i$, $f_i$ and $h_i$,
 $i=1,...,n$ be the
basic generators of ${\frak g}$ corresponding to the Cartan matrix.
 Let $(,)$ be the Killing form on ${\frak g}$.
Let ${\frak r}=\Sigma h_i\otimes h_i +\Sigma_{\Delta_+}
e_\alpha \otimes e_{-\alpha}+ \Sigma_{\Delta_-}
e_\alpha \otimes e_{-\alpha}$, where $\Delta_\pm$ are  the sets of all
the positive roots and the negative roots respectively,
 $e_{\alpha}$ is a root vector in ${\frak g}_{\alpha}$ and
$(h_i,h_j)=\delta_{ij}$, $(e_\alpha, e_{-\alpha})=1$.
${\frak r}$ is the  Casmir operator of ${\frak g}$. Let $\hat {\frak
g}$ be  the affine Lie algebra associated with ${\frak g}$. $
\hat {\frak g}$ has a concrete  realization as $\hat {\frak g}=
{\frak g}[x,x^{-1}]+\Bbb C[C]$, where $C$ is the central element.

\proclaim{Definition 3.5}Let $z$ be a formal variable. We define
the  elements  $\hat r$ and
 $r(z)$ in $\hat {\frak g}\otimes \hat {\frak g} $  and
 $\hat {\frak g}\otimes \hat {\frak g} [z,z^{-1}]$  as
 $$\hat r= \Sigma_{i,m} h_ix^m\otimes h_ix^{-m}
+\Sigma_{\alpha\in \Delta, i, m}
e_\alpha x^m \otimes e_{-\alpha}x^{-m},\tag 3.13   $$
$$r(z)=\Sigma_{i,m} h_ix^m\otimes h_ix^{-m}z^{-m}+
 \Sigma_{\alpha \in \Delta, i, m} e_\alpha x^m \otimes
e_{-\alpha}x^{-m}z^{-m}.$$
\endproclaim

We see that $\hat r$ is basically  like a
Casimir operator.

Let $e_i, f_i,  h_i$( $i=0,1,..,n$)
be the basic generators of $\hat {\frak g}$ for the corresponding
Cartan matrix of $\hat{\frak g}$.  Let $M$ be any finite dimensional
module of
$\hat {\frak g}$, Let $V_{\mu, k}$ be a highest weight module with
the  highest weight $\mu$ and
central extension $k$  of $\hat {\frak g}$. $k$ is a complex number.
$V_{\mu, k}$ is a graded module  such that $e_0$ and $f_0$ changes
the degree by $+1$ and $-1$ respectively as we explained before.

\proclaim{Theorem 3.2}
 $\pi_{V_{\mu,k}}\otimes \pi_M(r) $
maps $ V_{\mu,k}\otimes M$ to
$V_{\mu,k}\hat \otimes M$, commutes with
$e_i,f_i, h_i$, for  $i\neq0$, and
$$[\pi_{V_{\mu,k}}\otimes \pi_M (e_0),\pi_{V_{\mu,k}}\otimes \pi_M( r)]=
-2k (\text{id})\otimes (\pi_M (e_0)),\tag 3.14  $$
$$ [\pi_{V_{\mu,k}}\otimes \pi_M(f_0),\pi_{V_{\mu,k}}\otimes \pi_M
(r)]= 2k
(\text{id})\otimes  (\pi_M(f_0)). $$
By $V_{\mu,k}\hat \otimes M$, we mean the set of vectors in the form
of $\Sigma_{n\leq 0}  \mu(n)\otimes m_i $, where
$\mu(n)$ is a vector in $V_{\mu,k}$ of degree $n$.

\endproclaim

This follows from direct calculation, which is also a direct  corollary
of the corresponding assertion in the quantum case.

 Let $V=\Bbb C^{n} $ be the fundamental
representation of ${\frak g}$ as the case for
$\hat {\frak gl}(n)$ explained above. This representation can be extended to
 a representation of $\hat {\frak g}$.
It is clear that the concrete realization of
 $\pi_{V_{\mu,k}}\otimes \pi_V (r(z))$ can give us
explicitly the construction of the representation. That means,
for a specific representation $V_{\mu,k}$,  constructing  a representation
is equivalent to giving an  explicit expression of
$\pi_{V_{\mu,k}}\otimes \pi_V (r(z))$.  This is the central idea to
understand the classical spinor  constructions.
{}From now on in this section, we assume that ${\frak g}={\frak
gl}(n)$.  Let $V=\Bbb C^{n}$ be
a module of $\hat {\frak gl}(n)$ as we defined in the section above.
 Let $t$ be a real number such that $|t|$ is less than
1.

\proclaim{Theorem 3.3}Let ${\frak F}$ be the
standard map $V^* \otimes V$ to $\Bbb C$.
$(\Psi(z)\otimes I)
\Psi^* (zt)- 1/(1-t){\text id}\otimes
F$ and  $\lim_{t\rightarrow 1}(\Psi(z)\otimes I) \Psi^*
(zt)-1/(1-t) {\text id}\otimes F $ are well defined.
As a map from $V_{bf} \otimes V$ to $V_{bf} \otimes V[z,z^{-1}]$,
 $$-\bar r(z)= ({\text
id}\otimes I \otimes {\frak F})(\lim_{t\rightarrow 1}(\Psi(z)\otimes
I\otimes I) \Psi^* (zt)\otimes I-1/(1-t) {\text id}\otimes F\otimes I
)$$
$$=\pi _{V_{bf}}\otimes \pi _{ V}(\hat r(z)). \tag 3.15$$
As a map from $V_{bf} \otimes V$ to $V_{bf}
\hat \otimes V$,
$$-\bar r=\pi _{V_{bf}}\otimes \pi _{ V}
(D_z^{-1}\otimes 1)(\hat r(z)). \tag 3.16$$.
\endproclaim

By the limit above, we mean that we would take the limit for each
homogeneous component separately.

{\bf Proof.} It is clear that the first assertion
implies the second one.  The  limit we  take above is equivalent to the
normal ordering defined in \cite{FF}. A direct calculation shows that
$\bar r$ satisfies the property (3.14) of $r$ on the tensor module. We can
show by calculation that
the images  of the difference of the degree zero terms  of
the highest weight vectors is zero. However we know that the difference
between $\bar r$ and $r$ is an intertwiner. Thus the difference is zero.
Therefore  they are equal.

\documentstyle{amsppt}
\subheading{ 4.  Quantum Spinor
representation of
$U_q(\hat {{\frak gl}(n)})$ and
 Boson-Fermion Correspondence}

We  assume $n>2$ for the same reason explained in the
previous  section.

We will proceed to  construct the spinor
representation and  explain  the quantum boson-fermion correspondence
 for $U_q(\hat {{\frak gl}(n)})$. The degeneration of such a abstract
construction provides us  the classical
boson-fermion correspondence as in the section above.

\proclaim{Proposition 4.1}
$U_q(\hat {\frak gl}(n))$ is isomorphic to an algebra generated
by $U_q(\hat {\frak sl}(n))$, $g(n)$, $n\ne 0$ and K, such that $K$ is
central, $g(n)$ commute with $U_q(\hat {\frak sl}(n))$ and
$[g(l), g(m)]= \delta_{l, -m}mC$, where $C$ is the central element of $U_q(\hat
{\frak sl}(n))$\cite{DF}.
\endproclaim

{\bf Proof}.  From $U_q(\hat {\frak sl}(n))$, we can obtain $\bar
L^{\pm}(z)$, which satisfy (1.13a). There exist $g^\pm(z)=e^{ \Sigma b(\pm)
g(\mp m )z^{\pm m}}e^{\pm 1/2K}$,$m>0$, such that the operators  $g^\pm(z)\bar
L^{\pm}(z)$ stisfy the relation (2.10). This gives us an algebra
homomorphism. The isomorphism follows from the same proof as in
\cite{DF}.

With the proposition above and
 Lusztig's theorem about deformation of the
category of highest weight representations\cite{L}, we have that
the module  $V_{bf}$ can be deformed.
There is such  a   construction of this  module by
Frenkel and Jing \cite{FJ}.

As in the last section, we will start from the abstract construction.
We will denote this deformed  module of $V_{bf}$
 by $V_{\frak BF}$.
 $V_{\frak BF}$ as a module of ${U_q(\hat {\frak gl}(n))}$
can be decomposed
into  irreducible components: $V_{\frak BF}=\Sigma_{l\in \Bbb Z}
 \oplus V^l,$
where $V^l$ is   an  irreducible component  of  $V_{\frak BF}$, which
is a deformation of the module $V^l$ of
$\hat {\frak gl}(n)$ in section above.
But $V^l$ as a module of ${U_q(\hat {\frak sl}(n))}$ has
 infinite copies of $V_i$, which is a highest weight
representation of ${U_q(\hat {\frak sl}(n))}$ with the i-th
fundamental weight and central extension 1.
Let $V_{\frak bf}= \Sigma\oplus  V_i$. Surely, we assume that $V_i$
has the same grading shifting as in the classical case.

We now will present the evaluation  representations of $\hat{U_q({\frak
sl}(n))}$ on the fundamental representations of
the subalgebra $U_q({\frak sl}(n))$
 generated  $E_i, F_i, K_i$, $i\ne 0$ from  \cite{DO}.

We define a characteristic function $\theta^J(j)$ of a set $J$
 by $\theta^J(j)=1\,(j\in J),0\,(j\notin J)$. If $J$ is omitted, it
should be understood as $\Bbb Z_{\ge0}$.

Fix a positive integer $k$ such that $1\le k\le n-1$. Let $I=\{i_1,\cdots,i_k
\}$ be a subset of $\{0,1,\cdots,n-1\}$. For $I=\{i_1,\cdots,i_k\}$
we put $s(I)=i_1+\cdots+i_k$. Consider a
vector space $V^{(k)}$ spanned by the vectors $\{v_I\}$. We introduce an
algebra homomorphism $\pi^{(k)}:U_q(\hat {\frak sl}(n))\longrightarrow
{\text End}(
V^{(k)})$ by
$$ \pi^{(k)}(E_i)v_I=v_{I\backslash  \{i\}\cup \{i-1\}},$$
$$ \pi^{(k)}(F_i)v_I=v_{I\backslash \{i-1\}\cup  \{i\}}, \tag 4.1$$
$$\pi^{(k)}(K_i)v_I=q^{\theta^I(i-1)-\theta^I(i)}v_I, $$
for $i\neq0$, and
$$ \pi^{(k)}(E_0)v_I=v_{I\backslash \{0\} \cup \{n-1\}}, $$
$$\pi^{(k)}(F_0)v_I=v_{I\backslash \{n-1\}\cup \{0\}}, \tag 4.2 $$
$$\pi^{(k)}(K_0)v_I=q^{\theta^I(n-1)-\theta^I(0)}v_I. $$
Here $v_{I\backslash \{i\}\cup  \{i-1\}}$ should be understood as $0$ if
$i\notin I$ or
$i-1\in I$.
As a module $V^{(k)}$ of $U_q({\frak sl}(n))$ %
is isomorphic to the irreducible
highest weight module with highest weight corresponding to the
k-th node of the Dynkin diagram of ${\frak sl}(n)$.

Let  $V^{(k)}_z=V^{(k)}\otimes \Bbb C [z,z^{-1}]$. %
We can lift $\pi^{(k)}$ to %
an algebra homomorphism $\pi^{(k)}_z:U_q(\tilde {\frak sl}(n))
\longrightarrow {\text End}(
V^{(k)}_z)$ as follows:
$$\pi^{(k)}_z(E_i)(v_I\otimes z^n)=
\pi^{(k)}(e_i)v_I\otimes z^{n+\delta_{i0}},
\tag 4.3 $$
$$\pi^{(k)}_z(f_i)(v_I\otimes z^n)=\pi^{(k)}(f_i)v_I\otimes z^{n-\delta_{i0}},
$$
$$
\pi^{(k)}_z(t_i)(v_I\otimes z^n)=\pi^{(k)}(t_i)v_I\otimes z^n,$$
$$\pi^{(k)}_z(q^d)(v_I\otimes z^n)=v_I\otimes(qz)^n.$$
Then we have the following isomorphism of $U_q(\hat {\frak sl}(n))$--modules:
$$C_\pm^{(k)}:\quad V^{(k)}_{z(-q)^{\mp n}}
{\buildrel \sim\over\longrightarrow}
\Bigl({V^{(n-k)}_z}\Bigr)^{*a^{\pm1}}$$
$$ v_I \mapsto (-q)^{\pm s(I)}v_{I^c}^*.
\tag 4.4 $$
Here $I^c$ denotes the complement set of $I$ in $\{0,1,\cdots,n-1\}$ and
$\{v_{I^c}^*\}$ signifies the dual basis of $\{v_{I^c}\}\subset V^{(n-k)}$.
This gives us the evaluation representations we need.
By $(V^{(k)})^{*a^{\pm1}}$, we mean the right and the left dual of
$V^{(k)}$ respectively.

Let  $\Phi$ be an  intertwiner:
$$\Phi:V_{\frak bf} \longrightarrow V^1\otimes  V_{\frak bf},$$
$$\Phi (x)=e_1\otimes \Phi_1(x)+...+ e_n \otimes \Phi_n(x),\tag 4.5  $$
where $x\in V_{\frak bf}$ and $\{ e_i\}$ is a  basis for $V$.

Let $\Phi^*$ be the intertwiner:
$$\Phi^*:V_{\frak bf} \longrightarrow   V^{1*}
\otimes V_{\frak bf},$$
$$\Phi^* (x)=\bar e_1^*\otimes  \Phi_1^*(x)+...+\bar e_n^*\otimes
\Phi_n^*(x),\tag 4.6  $$
where $V^{1*}$ is the left  dual representation of $V$ of $\hat{U_q({\frak
sl}(n))}$, $x\in V_{\frak bf}$ and
 $\{ \bar e_i^*\}$ is a  basis for $V^{*1}$.

Let $\bar \Phi^*$ be the intertwiner:
$$\Phi^*:V_{\frak bf} \longrightarrow   ^*V^1\otimes  V_{\frak bf},$$
$$\bar \Phi^* (x)=  e_1^*\otimes \bar \Phi_1^*(x)+...+e_n^*\otimes
\bar \Phi_n^*(x),\tag 4.7  $$
where $^*V^1$ is the right  dual representation of $V$ of $\hat{U_q({\frak
sl}(n))}$, $x\in V_{\frak bf}$ and $\{  e_i^*\}$
 is a  basis for $^*V^1$.

{}From now on, we will identify  $V$ with  $V^1$.

There is an isomorphism from $^*V$ to $V^*$ by $a\rightarrow
q^{2\rho}a$, for $a\in ^*V$\cite{FR}, where
$\rho$ is the half  sum of all the positive roots of the $U_q({\frak
sl}(n))$ in $U_q(\hat {\frak sl}(n))$ generated by $E_i, F_i$  and$
K_i(i\neq0)$.

As in the previous  case, we can identify $V^*\otimes V$ with End($V)$.
With this identification, we define an operator $\tilde  L\in$
End$(V)$ $\otimes$ End$(V_{\frak bf})$:
$$\tilde {\frak L}= (\tilde {\frak L}_{ij})=
 ( (D^{2}_{q^{}}\Phi_j)\bar \Phi_i^*).
\tag 4.8 $$

$$\tilde {\frak L} \Delta (a)=
 ((  D^{2}_q\otimes 1 ) \Delta(a)) \tilde {\frak L} .\tag 4.9$$

The shift comes from the shift of  $\Phi_{q^{2}}$.

Here we have  a problem  that whether
the multiplication of two operators $D^{}_{q^{2}}\Phi_i$
 and $\bar \Phi_j^*$ is well defined, which we do not have in the finite
dimensional case. We procede to deal with this problem.
With  the condition  $|q|<1$, this
multiplication
is well defined, when we use Corollary  4.2 below, which comes from
 the results of  the correlation functions
of those intertwiners \cite{DO}.  We  assume, from now
on,  that
$|q|<1$.

Let $\Delta_j=j(n-j)/2n$, for $j=0,...,n-1$ and we extend this index cyclically
by $\Delta_j= \Delta_{j+n}$.

Let $\bar V_j$, $j=0,1,2,..,n-1$ equivalent to $V_j$ as a module of
$U_q(\hat {\frak sl}(n))$. We extend the index cyclically by
identifying $\bar V_j$ with  $\bar V_{j+n}$.
For an integer $l$,  let
$\bar {l}$ stand for the integer such that $l\equiv\bar{l}\mod n$,
$0\le\bar {l} <n$. Set
$$ I_{jk}=\{\overline  {j-k},\overline {j-k+1},\cdots,\overline {j}\}.$$

Let $\Phi_{\bar V_j}^{V^{(k)}\bar V_{j-k}}(z)$ denotes an intertwiner
from ${\bar V_j}$ to $V_z^{(k)}\otimes \bar V_{j-k}$.
Then our normalization reads as follows ($0\le j<n$):
$$\Phi_{\bar V_j}^{V^{(k)}\bar V_{j-k}}(z)|v_j\rangle
=z^{\Delta_{j-k}-\Delta_j}
v_{I_{jk}\setminus\{j\}}\otimes|v_{j-k}\rangle+\cdots,$$
$$\Phi_{\bar V_j}^{(V^{(k)})^{*a^{\pm}}\bar V_{j+k}}(z)|v_j\rangle
=z^{\Delta_{j+k}-\Delta_j}
v_{(I_{j,n-k}\setminus\{j\})^c}^*\otimes|v_{j+k}\rangle+\cdots, $$
where $v_i$ is  the highest weight vector
in $\bar V_i$.

Let $$ \Phi(z)=(1\otimes D^{-1}_z)\Phi=\Sigma
\oplus \Phi_{\bar V_j}^{V^{(1)}\bar V_{j-1}}(z)
, $$
$$  \Phi^*(z)= (1 \otimes D^{-1}_z)\bar \Phi^*=
\Sigma \oplus\Phi_{\bar V_j}^{(V^{(1)})^{*a^{+}}\bar V_{j+1}}(z),  $$
$$ \bar \Phi^*(z)= (1 \otimes D^{-1}_z)\Phi^*=
\Sigma \oplus\Phi_{\bar V_j}^{(V^{(1)})^{*a^{-}}\bar V_{j+1}}(z).\tag 4.10 $$
 We also have that
$$ \bar \Phi^*(z)=1\otimes q^{-2\rho} \Phi^*(zq^{-2n}).  \tag 4.11 $$

The  matrix coefficients of the highest weight vector
were obtained in \cite{DO}. The normalization above is given
in \cite{DO}.

We will present the commutation relations between those intertwiners.
In order to do this, we need to use  the
correlation functions given by Date and Okada \cite{DO}.

Let $$(z;p)_\infty=\prod_{j=0}^\infty(1-zp^j), $$

{}From now on, we will always use $\langle, \rangle$ to denote the
matrix coefficient of the corresponding highest weight vectors  of
the highest weight modules.

\proclaim{Proposition 4.2}\cite{DO} $$
\langle{\Phi_{\bar V_{j+k}}^{V_2^{(k)}\bar V_j}(z_2)
\Phi_{\bar V_{j+k+k'}}^{V_1^{(k')}\bar V_{j+k}}(z_1)}\rangle $$
$$=z_1^{\Delta_{j+k}-\Delta_{j+k+k'}}z_2^{\Delta_{j}-\Delta_{j+k}}
\prod_{i=1}^m{((-q)^{2i+|k-k'|-2}z_1/z_2;q^{2n})_\infty
\over((-q)^{-2i-|k-k'|}z_1/z_2;q^{2n})_\infty} \tag 4.12 $$
$$\qquad\quad\times\sum_{I\subset I_0^{(k+k')},|I|=k}
(-q)^{s(I_0^{(k)})-s(I)}(z_1/z_2)^{\mu_j(I)}
v_{(I_0^{(k+k')}\setminus I)[j]}\otimes v_{I[j]}.
$$
$$m=\min(k,k',n-k,n-k'), $$
$$\mu_j(I)=\sharp\{i\in I|i+j\ge n\}-\sharp\{i\in I_0^{(k)}|i+j\ge n\},$$
$$I[j]=\{\overline{i_1+j},\cdots,\overline{i_k+j}\}
\qquad\hbox{for }I=\{i_1,\cdots,i_k\}.$$
Here $I_0^{(k)} =\{0,1,\cdots,k-1\}$, if $k+k'>n$, we formally
understand $I_0^{(k+k')}=\{0,\cdots,
n-1\}\sqcup I_0^{(k+k'-n)}$, where the elements $0,\cdots,k+k'-n-1$ have
multiplicities, and we assume $I_0^{(k+k'-n)}\subset I,I_0^{(k+k')}\setminus
I$.
\endproclaim

Let
$$P_i:V^{(k)}\otimes V^{(k')}\longrightarrow V^{(k)}\otimes V^{(k')}
$$
be a linear map such that $P_i$ keeps the
copy of $V^{(i)}$ in $V^{(k)}\otimes V^{(k')}$ invariant
and map other irreducible components of
$U_q({\frak sl}(n))$ to zero.

We define matrices,  on
$V^{(k)}\otimes V^{(k')}$, a matrix  $\bar R_{V^{(k)}V^{(k')}}(z)$ as
$$P\bar R_{V^{(k)}V^{(k')}}(z)=\sum_{i=\max(0,k+k'-n)}^{\min(k,k')}\rho_i(z)
P_i,\tag 4.14
$$
where $P$ denotes the permutation.
The coefficients $\rho_i(z)$ satisfy
$$
{\rho_{i-1}(z)\over \rho_i(z)}={z-(-q)^{k+k'-2i+2}\over 1-z(-q)^{k+k'-2i+2}}.
$$
 $\rho_{\min(k,k')}(z)=1$ by our normalization of
$\bar R_{V^{(k)}V^{(k')}}(z)$.

Let
$$\Theta_p(z)=(z;p)_\infty(pz^{-1};p)_\infty(p;p)_\infty.$$

We define
$$R_{kk'}(z)=\rho^{(k,k')}(z)\bar R_{V^{(k)}V^{(k')}}(z),\tag 4.15 $$
$$ \rho^{(k,k')}(z)=z^{-kk'/n+\min(k,k')}
{((-q)^bz^{-1};q^{2n})_\infty((-q)^sz;q^{2n})_\infty\over
((-q)^bz;q^{2n})_\infty((-q)^sz^{-1};q^{2n})_\infty}
\prod_{i=1}^m{\Theta_{q^{2n}}((-q)^{2i+b}z^{-1})\over\Theta_{q^{2n}}((-q)^{2i+b}z)},
$$
where $b=|k-k'|,s=\min(k+k',2n-k-k')$,  $m$ is defined as
min$\{ k, k', n-k, n-k\}$  and $k,k'$=$1$ or $n-1$.

Let
$$
R_{kk'}^*(z)=({\text id}\otimes C_-^{(n-k')})R_{k,n-k'}(z(-q)^{-n})
({\text id}\otimes C_-^{(n-k')})^{-1},\tag 4.16 $$
$$ R_{kk'}^{**}(z)=(C_-^{(n-k)}\otimes C_-^{(n-k')})R_{n-k,n-k'}(z)
(C_-^{(n-k)}\otimes C_-^{(n-k')})^{-1},$$.

\proclaim{Corollary  4.2} Let $k,k'$=$1$ or $n-1$. Then
$$\langle{\Phi_{\bar V_{j+k}}^{V_2^{(k)}\bar V_j}(z_1)
\Phi_{\bar V_{j+k+k'}}^{V_1^{(k')}\bar V_{j+k}}(z_2)}\rangle
=PR_{kk'}(z_1/z_2)\langle{
\Phi_{\bar V_{j+k'}}^{V_2^{(k')}\bar V_j}(z_2)
\Phi_{\bar V_{j+k+k'}}^{V_1^{(k)}\bar V_{j+k'}}(z_1)}\rangle,
$$
$$\langle {
\Phi_{\bar V_{j+k}}^{V_2^{(k)}\bar V_j}(z_1)
\Phi_{\bar V_{j+k-k'}}^{V_1^{(k')*a^{-1}}\bar V_{j+k}}(z_2)}\rangle $$
$$
\qquad\qquad=PR_{kk'}^*(z_1/z_2)\langle{
\Phi_{\bar V_{j-k'}}^{V_2^{(k')*a^{-1}}\bar V_j}(z_2)
\Phi_{\bar V_{j+k-k'}}^{V_1^{(k)}\bar V_{j-k'}}(z_1)}\rangle, $$
$$\langle{\Phi_{\bar V_{j-k}}^{V_2^{(k)*a^{-1}}\bar V_j}(z_1)
\Phi_{\bar V_{j-k-k'}}^{V_1^{(k')*a^{-1}}\bar V_{j-k}}(z_2)}\rangle
\tag 4.17 $$
$$
\qquad\qquad=PR_{kk'}^{**}(z_1/z_2)\langle{
\Phi_{\bar V_{j-k'}}^{V_2^{(k')*a^{-1}}\bar V_j}(z_2)
\Phi_{\bar V_{j-k-k'}}^{V_1^{(k)*a^{-1}}\bar V_{j-k'}}(z_1)}\rangle. $$

In the neighborhood of $|z_1/z_2|=1$, both sides of the second formula
above  with $k=k'$ have a simple pole at $z_1=z_2$. Its residue is given by
$$P{\text Res}_{z_1=z_2}\langle{
\Phi_{\bar V_{j+k}}^{V_2^{(k)}\bar V_j}(z_1)
\Phi_{\bar V_j}^{V_1^{(k)*a^{-1}}\bar V_{j+k}}(z_2)\,d\left({z_1\over
z_2}\right)
}\rangle $$
$$=h^{(k)}\sum_{I\subset I_0^{(n)},|I|=k}v_I\otimes v_I^*,\tag 4.18 $$
where
$$h^{(k)}={(q^{2n-2};q^{2n})_\infty\over(q^{2n};q^{2n})_\infty}
\prod_{i=1}^{\min(k,n-k)-1}
{(q^{2n-2i-2};q^{2n})_\infty\over(q^{2i};q^{2n})_\infty}.
$$
\endproclaim

Let

    $$f(z_1/z_2)= {(z_2/z_1;q^{2n})_\infty
\over(q^{-2}z_2/z_1;q^{2n})_\infty},$$
    $$F(z_1/z_2)={((q)^{2n-2}z_2/z_1;q^{2n})_\infty
\over(q^{2n} z_2/z_1;q^{2n})_\infty}. $$
Note that $(1-z_2/z_1q^{-2})f(z_1/z_2)=(1-z_2/z_1)/F(z_1/z_2)$

Let $R_{1,1}(z_1/z_2)$, $R^{**}_{n-1,n-1}(z_1/z_2)$ and
$R^{*}_{1,n-1}(z_1/z_2)$  be  matrices as  defined above
but the first term of $\rho^{(k,k')}(z)$ in (4.15) will be substitued by
$z^{\delta_{k,k'}}$.

\proclaim{Corollary 4.2}\cite{FR}\cite{ DFJMN} \cite{DO}
Let $z$,  $z_1$ and $z_2$ be  formal variables.
 $\Phi(z)$ and $ \bar \Phi^*(z)$ satisfy
the commutation relations:
$$1/f(z_1/z_2)\Sigma\Phi_j(z_1) \Phi_i(z_2)e_i\otimes e_j=$$
$$P'
1/f(z_1/z_2)(R_{1,1}(z_1/z_2))( \Sigma \Phi_{j'}(z_2)
\Phi_{i'}(z_1)e_{i'}\otimes e_{j'}), $$
$$1/f(z_1/z_2)\Sigma \bar \Phi^*_i(z_1)\bar \Phi^*_j(z_2)
e_j^*\otimes
e_i^*=$$
$$ P'
1/f(z_1/z_2)(R^{**}_{n-1,n-1}(z_1/z_2)
 (\Sigma
\bar \Phi^*_{i}(z_2)\bar \Phi^*_{j}(z_1)e_j^*\otimes
e_i^*), \tag 4.19 $$
$$1/F(z_1/z_2)\Sigma \bar \Phi_j(z_1)\Phi^*_i(z_2)e^*_i\otimes
e_j-1/(1-z_1/z_2)F=$$
$$
P'1/F(z_1/z_2)(R^{*}_{1,n-1}(z_1/z_2))
 (\Sigma\Phi_{j}^*(z_1) \bar \Phi_{i}(z_2)
e_i\otimes e^*_j)-(z_1/z_2)/(z_1/z_2-1)F,
$$
$$\lim_{z_1\rightarrow 1}\lim_{z_1\rightarrow z_2, |z_1|<|z_2|}(z_1-z_2)
 \Sigma \Phi_i(z_1)\bar \Phi^*_i(z_2) e_i\otimes e_i^*=
P {(q^{2n-2};q^{2n})_\infty\over(q^{2n};q^{2n})_\infty}F,$$
where  $f(z_1/z_2)$ and $F(z_1/z_2)$ are expanded in
the power series of $z_2/z_1$ on the left hand side of the
the formulas above but in the power series of $z_1/z_2$ on the
right hand side.  $F=\Sigma e_i\otimes  e^*_i$.
\endproclaim

{\bf Proof.}Our arguement is  based on the formlas in \cite{DO} of the
correlation functions.

The argument for the
first two formulas  is straightforward.
First, we know that the formula is true on the
level  of correlation functions  of the highest weight vectors of
$V_{\frak bf}$,
 due to the fact that after we factor  out
those functions  as above, the correlation functions  of
the operators on both sides of
highest weight  vectors  are polymonials as showed in \cite{DO}.
On the other hand,  both sides are intertwiners, this relation thus can be
proved to be true for the matrix coefficient of any
two vectors,  thus both sides are equal.

 As for the formula for the commutation relations between $\Phi_i$ and
$\bar \Phi_j^*$, the argument goes as follows:
the first part is the same as that of above, namely,
the formula is valid on the level  of correlation functions
of the highest weight vectors;
secondly,  $\Sigma_{n\in \Bbb Z} (z_1/z_2)^nF $ is an
invariant vector, thus ${\text id}\otimes \Sigma_{n\in \Bbb Z} (z_1/z_2)^nF $
is also an intertwiner. Thus the difference of
two sides is also an intertwiner, then we can show that
the third formula  is valid on the level of correlation functions
of any two vectors of the two sides of the third formula.
Thus it is valid.

 Basically the idea appeared  in \cite{FR}
 \cite{ DFJMN}, which is to
study  the relations between the correlation functions.

The commutation relation of each homogeneous component
of $\Phi_i(z)$,  $\bar \Phi_j^*(z)$ described by using $R$-matrix
will degenerate into the commutation relation  with $\delta (z)$.

 The locations  of
the poles of the
correlation functions of $\Phi_j(z_1)\bar \Phi_i^*(z_2) $
do not include the line  $z_1q^2=z_2$. From the commutation
relations and the condition that  $|q|<1$,
the multiplication of $\bar \Phi^*(z)$ and
$\Phi(zq^2)$  is well  defined.
Thus  $(D^{}_{q^{2}}\Phi_i)\bar \Phi^*_j e^*_j\otimes e_i =
(1\otimes 1\otimes D_z^{-1})
(1\otimes \Phi(zq^2))\bar \Phi^*(z)$ is well defined.

We know that both  $V \otimes  V_i$ and  $V_{q^{2}}\otimes V_i$
  are irreducible
\cite{FR} \cite{KKMMNN} \cite{JM}.
This shows that the dimension of the space of the operators
 $X: V\otimes  V_{\frak bf} \longrightarrow  V_{q^{2}}\otimes
  V_{\frak bf},$
which  satisfy the following relation :
$$ X \Delta(a)= ( D^{2}_q\otimes 1 ) \Delta (a) X, $$
is $n$. Thus the difference between $\tilde L$ and
$L$
is  a constant factor if we restrict it to each irreducible component
of $V\otimes V_i$, which can  be determined  by looking at their
actions  on the highest weight vectors. By looking at
the  homogeneous component of
degree 0 of the correlation functions,
 we know that there is   a universal
factor $c$ for all the irreducible components, which can be derived
by comparing the actions  of $\frak L$ and $\tilde L$ on the
highest weight  vectors.

Let $\tilde {\frak L}(z)$=$(1\otimes \Phi(zq^2))\bar \Phi^*(z)$.

\proclaim{Theorem 4.1}
 $$  {\frak L}(z)=c (D_{z}\otimes 1) \tilde {\frak  L},   \tag 4.20  $$
where $c=\frac {{\text tr}(v_0,(D^{}_{q^{2}}\bar \Phi^*) \Phi v_0)} {
  {\text tr}(v_0, {\frak L}(1)v_0)} $ and  $v_0$ is the highest weight
vector of $V_0$.
\endproclaim

We can remove the $c$ in (4.11)
 by renormalizations  of $\Phi$ and $\bar \Phi^*$.

{}From the definition of ${\frak L}(z)$, we know
that  $\tilde {\frak L}(z) $ can be decomposed into the
product of ${\frak L}^+(z)$ and ${\frak L}^{-}(z)^{-1}$,  which
reminds us the polar  decomposition. With the
fermionic realizations, we  could factor out
${\frak L}^+(z)$
easily by looking at $\tilde {\frak L}(z) v_i$, where $v_i$ are the highest
weight  vectors for
$V_{\frak bf}$. This decomposition is unique.

Now,    let's consider
$\hat{U_q({\frak gl}(n))}$. As in \cite{DF},  by adding an extra Heisenberg
algebra $h(n)$, which commutes
 with  $U_q(\hat {\frak
sl}(n))$, we would obtain the representation of
$\hat{U_q({\frak gl}(n))}$. From \cite{DF} and (1.11b)
\proclaim{Proposition 4.2}There exists complex numbers
$a(n)$ such  that $$\bar { {\frak  L}}^+(z)=({{\frak L}^+}(z)
\otimes e^{\sum_{m \in \Bbb Z_{+}}
a(m)h(-m)(z)^{m}})^{-1}, $$
 $$\bar { {\frak  L}}^{-}(z)=( {\frak L^{-}}(z)\otimes
e^{ \sum_{m\in \Bbb Z_{+}}a(-m)h(m)z^{-m}})^{-1}, \tag 4.21  $$
satisfy the commutation relation (2.10).
\endproclaim

These operators acting on the tensor of the Fock space of
$H(-m)$ and $V\otimes V_{\frak bf}$.
We will denote it by  $V_{gl}$. $V_{gl}=V_{\frak bf}\otimes  V_{h} $,
where $V_h$ is the module generated by the extra Heisenberg algebra
of  ${U_q(\hat {\frak gl}(n))}$. Let's give the same Gauss
decomposition to these new operators. Then the action of the product
of the zero component of the diagonal components  of their
decomposition, which we denote by $T$ and $T^{-1}$, are  $1$.
Let $A$ be a group algebra generated by
a lattice $\Bbb Z a$.  Let $\bar V=\Sigma \oplus  V^i_{gl}\otimes e^{mn+i}$,
$m\in
\Bbb Z$ in the space $V^i_{gl}\otimes A$, where
$V^i_{gl}=V_i\otimes V_h$. We define
$\bar V$ to be a module of $\hat{U_q({\frak gl}(n))}$, such that
all other elements  acting  only on $V_{gl}$, but the action of
$T$ on $ V_{gl}\otimes e^{ma}$ is a multiplication of $q^m$.

\proclaim{Proposition 4.4} $T^{\pm1/2}\bar {\frak L}^\pm(z)$ gives us
 representation $L^\pm(z)$ of $\hat{U_q({\frak gl}(n))}$ on
 $\bar V$. $\bar V$  is equivalent to
$V_{\frak BF}$.
\endproclaim

Let $\bar h(z)=\Sigma a(-m) h(m)z^{-m}/(q^2-1)$ and $\bar h(z_1)\bar
h(z_2)$= $g(z_1/z_2)+:\bar h(z_1)\bar  h(z_2):$
Let $G(z_1/z_2)=e^{g(z_1/z_2)}$. Then $:e^{\bar h(z_2)}::e^{\bar h(z_1)}:$
$=G(z_1/z_2):e^{\bar h(z_2)+\bar h(z_1)}:$.

\proclaim{Definition 4.1}Let
$z$,  $z_1$ and $z_2$ be formal variables.
  Affine  Quantum Clifford algebra is  defined as
 an  associative algebra generated by
$\psi(z)=(\psi_i(z))=(\Sigma \psi_i(m)z^{-m})$ and $
\psi^*(z)=(\psi^*_i(z))
=(\Sigma \psi^*_i(m)z^{-m})$, $0<i<n+1$
 satisfying
the commutation relations:
$$1/(f(z_1/z_2)/G(z_1/z_2))\Sigma\psi_j(z_1) \psi_i(z_2)e_i\otimes e_j=$$
$$
P'1/(f(z_1/z_2)/G(z_2/z_1))(R_{1,1}(z_1/z_2))( \Sigma
\psi_{i'}(z_2)\psi_{j'}(z_1)e_{j'}\otimes e_{i'}),$$
$$1/(f(z_1/z_2)/G(z_1/z_2))\Sigma \psi^*_j(z_1)\psi^*_i(z_2)e_i^*\otimes
e_j^*=$$
$$P'1/(f(z_1/z_2)/G(z_2/z_1))(R^{**}_{n-1,n-1}(z_1/z_2)
 (\Sigma\psi^*_{i}(z_2)
\psi^*_{j}(z_1)e_j^*\otimes e_i^*), \tag 4.22  $$
$$1/(F(z_1/z_2)G(z_1/z_2))\Sigma \psi_j(z_1)\psi^*_i(z_2)e^*_i\otimes \bar
e_j-1/(1-z_2/z_1)F=$$
$$P'1/(F(z_1/z_2)G(z_2/z_1))(R^*_{1,n-1}(z_1/z_2))
 (\Sigma \psi^*_{j}(z_1)\psi_{i}(z_1)\bar
e_i\otimes e^*_j)$$
$$-(z_1/z_2)/(z_1/z_2-1)F,
$$
$$\lim_{z_1\rightarrow 1}\lim_{z_1\rightarrow z_2, |z_1|<|z_2|}1/(z_1-z_2)
 \Sigma \psi_i(z_1)\psi^*_i(z_2)e^*_i\otimes e_i =
 {(q^{2n-2};q^{2n})_\infty\over(q^{2n};q^{2n})_\infty}F/G(1),$$
where
the functions of  the left and the right  hand sides
 are expanded
in $z_2/z_1$ and $z_1/z_2$ respectively.
\endproclaim

\proclaim{Theorem 4.2}
Quantum Clifford algebra  is isomorphic to the algebra generated by
$\Phi(z)\otimes e^{-\bar h(z)}\otimes e^{-a}$ and
$\Phi(z)^* \otimes e^{\bar h(z)} \otimes e^{a}$ on $\bar V$.
\endproclaim

{\bf Proof.}  It is straight forward to show that the map from
$\psi(z)$ to $\Phi(z)\otimes e^{-\bar h(z)}\otimes e^{-a}$ and $\psi^*(z)$
to $\Phi(z)^* \otimes e^{\bar h(z)} \otimes e^{a}$ is a surjective
algebra homomorphism. Because  $\rho^{kk}(z)$ for $k=1$ and $k=n-1$
has a factor in the form of $z (1-z^{-1})/(1-z)=-1$ and (2.9),
if we shift the degree of $\psi(z)$ and $\psi^*(z)$ by $\pm 1/2$ as
in the classical case, we can define the Fock space as in
Definition 3.3, then derive the character with the calculation
based on the R-matix on $\Bbb C^n\otimes \Bbb C^n$ given in Section 2
on the specific basis we chose.  By
comparing the character,  we can prove  the isomorphism.

In definition 4.1, we expect that the complicated
functions showing in Definition 4.1 will conceal each other, such that
we would get  manageable and easy formulas. If this hypothesis
is true, it hints that,  as in the classical case, we should look at
the
corresponding case of ${\frak gl}(n)$ instead of ${\frak sl}(n)$,
which should make the case  much simpler.

With Theorem 4.2, we actually can start from the abstract algebra
defined
in Definition 4.1. Then we can derive $L$, which leads to
the realization of $L^{\pm}(z)$ as  defined in Section 2.
{}From Theorem 2.1 in Section 2,    through the Gauss  decomposition of
${ L}^{\pm}(z)$,  we obtain all the quantum bosons out of
 ${ L}^{\pm}(z)$. Thus
we obtain  the realization of  the quantum boson-fermion
correspondence in one direction.
But we can not write explicit simple   formulas
due to the difficulties  coming from the polar
decomposition of ${L}(z)$ and the Gauss decomposition of
${L}^\pm(z)$, however the case for $U_q(\hat {\frak sl(2)})$
can be solved in a relatively easy way through computation.

 On the other hand, based on the
work Koyayma \cite {K}  and Frenkel-Jing construction,
we can write down partially  the realization the
quantum fermions in Bosons. With our results, we expect that
a complete formula is very possible, if we consider ${\frak gl}(n)$
instead of ${\frak sl}(n)$ as we explain in Remark 1.

Let $\bar F(z_1/z_2)= (1-z_2/z_1 q^{-2}) F(z_1/z_2)$
Let $H(z) = \Sigma H(n)z^{-n}$, $n\neq 0$ be an Heisenberg algebra
such that $:e^{H(z_1)}: :e^{-H(z_2)}:= 1/(\bar
F(z_1/z_2)):e^{H(z_1)-H(z_2)}:$.
Let $\bar \Psi(z)=\Phi(z)\otimes :e^{-H(z)}:$, $\bar
\Psi^*(z)=\Phi^*(z)\otimes :e^{H(z)}:$.  Let $\tilde V$ be the pace of
tensor of $V_{\frak bf}\otimes H$, where $H$ is the space generated by
$H(n),  n<0$. Then, on $\tilde V$, we have

\proclaim{Theorem 4.3}
Let $z$,  $z_1$ and $z_2$ be formal variables.
$\bar \Psi(z)$ and $\bar \Psi^*(z)$ satisfies the following relations:
$$\Sigma \bar \Psi_j(z_1) \bar \Psi_i(z_2)e_i\otimes e_j=$$
$$(1-z_1/z_2)
 P'\bar F(z_2/z_1)/(f(z_1/z_2))(R_{1,1}(z_1/z_2))( \Sigma
\bar \Psi_{i'}(z_2)\bar \Psi_{j'}(z_1)e_{j'}\otimes e_{i'}),$$
$$\Sigma \bar \Psi^*_j(z_1)\bar \Psi^*_i(z_2)\bar e_i^*\otimes
\bar e_j^*=$$
$$(1-z_1/z_2)P'\bar F(z_2/z_1)/(f(z_1/z_2))(R^{**}_{n-1,n-1}(z_1/z_2)
 (\Sigma \bar \Psi^*_{i}(z_2)
\bar \Psi^*_{j}(z_1) \bar e_j^*\otimes \bar e_i^*),
\tag 4.23
$$
$$\Sigma \bar \Psi_j(z_1)\bar \Psi^*_i(z_2)e^*_i\otimes \bar
e_j-1/(1-z_2/z_1)F(1-q^{-2})=$$
$$P'1/(F(z_1/z_2)\bar F(z_2/z_1)))(R^*_{1,n-1}(z_1/z_2))
 (\Sigma \bar  \Psi^*_{i}(z_2) \bar \Psi_{j}(z_1)\bar
e_j\otimes e^*_i)$$
 $$ -(z_1/z_2)/(z_1/z_2-1)(1-q^{-2})F,
$$
$$\lim_{z_1\rightarrow 1}\lim_{z_1\rightarrow z_2, |z_1|<|z_2|}1/(z_1-z_2)
 \Sigma \bar \Psi_i(z_1) \bar \Psi_i^*(z_2) e^*_i\otimes \bar e_i =$$
$$(1-q^{-2})
F,
$$
where
the functions of  the left and the right  hand sides
 are expanded
in $z_2/z_1$ and $z_1/z_2$ respectively.
\endproclaim

The proof comes from straight calculation.

  \documentstyle{amsppt}
{\bf Discussion.}

Theorem 4.3 gives a realization of an algebra, which has the same definiton
formulas for form factors
in quantum field theories\cite{Sm}, but the R-matrix here is different
with a function factor.

To construct local operators in the theory of fromal factors is
a very important problem\cite{Sm}. If we consider the case of
$U_q(\hat{\frak sl}(n))$,
and if we consider the intertwiners in
(1.18) and (1.19), which we will call right intertwiners,
 as the basic generators to define
form factors, it is known that certain compositon of
the inetrtwiners of type as in (1.14) and (1.15), which we
will call left intertwiners, gives
local operators. However, how to derive the left intertwiners
from the right intertwiners is a problem.
With the spinor constrution we derive above, it is clear that we can
derive the operator ${\frak L}^\pm(z)$ through Gauss decomposition of
${\frak L}(z)$ constructed out of the right intertwiners.
The proper composition of ${\frak L}^\pm(z)$ with the rihgt intertwiners
obviously gives us the left intertwiners. Thus our construction
provides a way to obtain local operators directly from the
algebra, which defines form factors.
Moreover, we propablely can
 modify the algebra with an extra Heisenberg algebra
similar to the case in Theorem 4.3, such that we can
derive in the same way operators $l^\pm(z)$ comimg from the Gauss
decomposition of the operator $l(z)$ built out of those
modified right intertwiners
to derive operators similar to the left intertwiners,
 which, however, simply  commutes with those modified right
intertwiners, which is basically the definiton property of
local operators.

As a natural continuation, we are expecting to  apply the
same idea to undercovering  the underlying structure of the
so-called vertex operator algebras, which
should lead us to the corresponding deformed structure,
an axiomatic formulation  quantum vertex operator algebra
via the  representation theory   and the  structure theory of
quantum affine algebras.
The complete establishment of such a theory should provide a
proper mathematical setting to understand the
massive quantum field theory in theoretical physics.

Recently, there appeared two papers by  E.  Stern \cite{St}
and M. Kashiwara,   T.  Miwa, and E.  Stern(q-alg/9508006),
where they study in detail basically the same Fock space as the Fock
space defined by the quantum Clifford algebra  in
this paper.

{\bf Acknowledgments}
This paper is part of a dissertation under Professor I. Frenkel
submitted to Yale University in May 1995. I would like to
 thank my advisor, Igor B.  Frenkel, for his
guidance and   his constant and creative  encouragement.
I would  like to thank Prof. M. Jimbo for his stimulating
discussion and advice, especially as  concerns
the commutation relations  of the
intertwiners. I would laso like to thank  the Sloane  foundation for
the dissertation fellowship.

\documentstyle{amsppt}

\end